\documentclass[12pt,a4paper]{article}  

\usepackage{amsmath,amsthm,amsfonts,amssymb,oldgerm}
\usepackage[frenchb,english]{babel}
\usepackage{indentfirst}
\numberwithin{equation}{section}
\usepackage[citecolor=blue,colorlinks=true]{hyperref}
\usepackage{color}
\allowdisplaybreaks[1]
\usepackage[top=3cm, bottom=3cm, right=2.2cm, left=3cm]{geometry}

\newcommand{\disp}{\displaystyle}

\newcommand{\E}{\mathds{E}}
\renewcommand{\H}{\mathds{H}}
\newcommand{\N}{\mathds{N}}
\renewcommand{\P}{\mathds{P}}

\newcommand{\R}{\mathds{R}}

\newcommand{\1}{\mathds{1}}

\newcommand\cA{{\mathcal A}}

\newcommand\cD{{\mathcal D}}
\newcommand\cE{{\mathcal E}}

\newcommand\fB{{\mathfrak B}}

\newcommand\fF{{\mathfrak F}}

\newcommand\fR{{\mathfrak R}}

\usepackage{dsfont}				
\usepackage{cleveref}

\newtheorem{theorem}{Theorem}[section]
\crefname{theorem}{Theorem}{Theorems}

\crefname{proof}{Proof}{Proofs}

\newtheorem{proposition}{Proposition}[section]
\crefname{proposition}{Proposition}{Propositions}

\newtheorem{lemma}{Lemma}[section]
\crefname{lemma}{Lemma}{Lemmas}

\crefname{equation}{equation}{equations}

\crefname{section}{Section}{Sections}

\newtheorem{remark}{Remark}[section]
\crefname{remark}{Remark}{Remarks}

\newtheorem{assumption}{Assumption}[section]
\crefname{assumption}{Assumption}{Assumptions}

\crefname{hypothesis}{Hypothesis}{Hypothesis}

\crefname{example}{Example}{Example}

\crefname{note}{Note}{Notes}

\newtheorem{corollary}{Corollary}[section]
\crefname{corollary}{Corollary}{Corollary}

\renewcommand{\abstractname}

\begin{document}

\begin{center}
\large \textbf{A PIECEWISE DETERMINISTIC LIMIT  FOR A MULTISCALE STOCHASTIC SPATIAL GENE NETWORK}
\end{center}

\begin{center}
\textsc{By Arnaud Debussche\footnotemark[1]  and  Mac Jugal Nguepedja Nankep}\footnotemark[1] \footnotetext[1]{Univ Rennes, CNRS, IRMAR - UMR 6625, F-35000 Rennes, France}
\end{center}




\begin{abstract}
We consider multiscale stochastic spatial gene networks involving chemical reactions and diffusions. The model is Markovian and the transitions are driven by Poisson random clocks. We consider a case where there are two different spatial scales: a microscopic one with fast dynamic and a macroscopic one with slow dynamic. At the microscopic level, the species are abundant and for the large population limit a partial differential equation (PDE) is obtained. On the contrary at the macroscopic level, the species are not abundant and their dynamic {remains governed} by jump processes. It results that the  PDE governing the fast dynamic contains coefficients which randomly change.  The global weak limit is an infinite dimensional continuous piecewise deterministic Markov process (PDMP). Also, we prove convergence in the supremum norm.
\end{abstract}

\section{Introduction}

Within the last decades, spatial and stochastic modeling has been widely used for systems of biochemical reactions. Commonly used models  describe systems where reactants undergo chemical reactions and can diffuse in the considered spatial domain. They are used either deterministically or stochastically. 

Deterministic models are reaction-diffusion equations, which are partial differential equations (PDEs). These equations are solved analytically, and/or simulated numerically in which case results are obtained relatively fast. 
However, these models capture a macroscopic dynamic and are valid only in high concentrations contexts. When some of the interacting species are present in small quantity,  a stochastic description seems more accurate but the direct computation of stochastic models is extremely time consuming. This problem is typical for  multiscale systems. A remedy is to compute abundant species as continuous variables that follows deterministic motions. Species in low number remain stochastic and are directly simulated, one talks about hybrid approximation.
Lately, hybrid algorithms have been proposed for the simulation of multiscale spatial models 
arising from cellular biology or related fields. We refer for instance to
\cite{DunErbZyg2016} for cell regulatory networks and enzyme cascades, or \cite{NoelKarCheSch2015} for molecular communication.

Stochastic hybrid systems, and especially piecewise deterministic Markov processes \- (PDMPs), form a class of systems that has get very popular throughout the past decade, as it proposes a quite natural simplification of multiscale systems. 

On the mathematical side, the theory of PDMPs was initiated in finite dimension by Davis \cite{Davis1984, Davis1993}.  Finite dimensional PDMPs are suited to model spatially homogeneous situations. Recently, Buckwar and Riedler \cite{Riedler2011bis, Riedler2011} have extended PDMPs to the infinite dimensional case, in order to model the propagation of action potentials in neurons (see also \cite{Riedler2012, Genadot2013}). Mathematical results about stochastic hybrid systems concern existence, well-definiteness, Markov property \dots

On another hand, very few results exist on the mathematical validity of these hybrid models.  Do they really approximate correctly or not multiscale stochastic sytems ?


In many situations,  laws of large numbers hold for one scale systems  - {homogeneous or with spatial dependence}. 
Among the numerous existing references, let us quote \cite{Kurtz1970} pioneering work in homogeneous framework, \cite{Arnold1980}, 
\cite{blount1992} in the context of chemical reactions with spatial diffusion modeling, or more recently \cite{NziPardouxYeo2019} for a compartmental SIR epidemic model.
These allow  to approach one scale stochastic systems, under the assumption of a large size of individuals, by a corresponding deterministic version.
Central limit theorems results have also been proved, as well as large or moderate deviations results. 
See, among others, \cite{Kotelenez1986} - for chemical reaction-diffusion models - with respect to the former group of results, and \cite{Yeo2019} - for epidemic models - with respect to the latter group of results.

Radulescu, Muller and Crudu in \cite{RadMulCru2007}, then Crudu, Debussche and Radulescu in \cite{Arnaud2009} have proposed hybrid approximations by finite dimensional PDMPs, for some multiscale stochastic homogeneous gene networks. Then a rigorous justification has been given in \cite{Arnaud2012}. 
The question is more complicated in the spatially inhomogeneous case, where there are much more possibilities of modeling,
as can be seen in a previous work (\cite{DebusscheNankep2017}), and in \cite{Nankep2018}, chapter 4 - or equivalently - \cite{Kouegou2019}, chapter 3. The authors of the two latter references have proved a law of large numbers and large deviation inequalities in a common work on a spatial model of cholera epidemic.
\par
In the previous work \cite{DebusscheNankep2017}, a multiscale system with spatial dependance was considered. As usual, the spatial domains is divided in a finite number of cells and in each cell two species are present but one - the continuous - is much more abundant. However one - the discrete - has only few individuals in each cells. 
It is proved in this 
work that the limit is a reaction-diffusion PDE for the continuous component modeling continuous species, coupled to an ordinary differential equation (ODE) driving the discrete component modeling species in low numbers. This may seem surprising and one might expect a PDMP at the limit {but this is in fact natural} 
because, even though the discrete specie is less abundant, its global number is of the order of the number of cells and grows to infinity and stochastic effects 
disappear.

Another situation is considered in the present paper. We consider the case when the discrete species have also a different spatial scale: their size is macroscopic. This may correspond to cells 
or to group of cells, depending on the context of application.
Then their number remains small and stochasticity remains at the limit. The limit is a continuous infinite dimensional PDMP whose continuous component satisfies, between the jumps of its discrete component, a reaction-diffusion equation parametrized by the value of the discrete component between the considered jumps, see section \ref{s2.3}  for the precise derivation of the stochastic model and section \ref{s2.4.1} for a formal derivation of
the limit model. 

The rest of this article is organized as follows. In section \ref{s2}, we briefly recall the definition of an infinite dimensional PDMP and collect useful results about it. Then, we develop our model of interest and present heuristics allowing to identify its limit. The main result of convergence is stated and proved in section \ref{s3}. The tightness of the process is proved first thanks to similar arguments as in \cite{DebusscheNankep2017},   then the limit is rigorously identified through the martingale problem, and we conclude by a truncation argument. This step is more difficult than in \cite{DebusscheNankep2017}, new difficulties appear. 
Section \ref{s4} contains the proof of uniqueness of the solution of a martingale problem associated to infinite dimensional PDMPs, a crucial result used in 
section \ref{s3}.

\begin{subsection}*{Acknowledgment}

A. Debussche and M. J. Nguepedja Nankep are partially supported by the French government thanks to the "Investissements d'Avenir"
program ANR-11-LABX-0020-01. 
\end{subsection}

\section{Modeling and asymptotics}\label{s2}

\subsection{Notation}

Let $(Z,\Vert\cdot\Vert_Z)$ and $(\tilde{Z},\Vert\cdot\Vert_{\tilde{Z}}$ be Banach spaces. The product space $Z\times\tilde{Z}$ is equipped with the norm $\Vert\cdot\Vert_Z+\Vert\cdot\Vert_{\tilde{Z}}$. 
The space of continuous linear maps from $Z$ to $\tilde{Z}$ is denoted by $\mathcal{L}(Z,\tilde{Z})$. If $Z=\tilde{Z}$, one simply writes $\mathcal{L}(Z)$.
The operator norm is denoted $\Vert\cdot\Vert_{Z\rightarrow\tilde Z}$ and when there is no risk of confusion, we denote it $\Vert\cdot\Vert$.
For $\tilde Z=\R$, $Z'$: the space of continuous linear forms on $Z$, {\it i.e.} the topological dual of $Z$.

$\fB(Z)$ (resp. $\fB_b(Z)$) is the space of Borel-measurable (resp. bounded Borel-measurable) real valued functions on $Z$. The space $\fB_b(Z)$ is endowed with the supremum norm \[ \Vert f\Vert_{\fB_b(Z)}=\sup_{x\in Z}|f(x)|=\Vert f\Vert_\infty. \]

$\disp C_b^k(Z)$, $k\in \N$ is the space of real valued functions of class $C^k$, i.e. $k$-continuously Fr\'echet differentiable,  on $Z$ which are bounded and have uniformly bounded succesive differentials. It is equipped with the norm \[ \Vert f\Vert_{C_b^k(Z)}=\sum_{i=0}^k\big\Vert D^{i}f\big\Vert_\infty, \] where $D^{i}f$ is the $i$-th differential of $f\in C_b^k(Z)$, and $\disp C_b^0(Z)=C_b(Z)$ is the set of bounded continuous real valued functions on $Z$. \par 

$C^{l,0}(Z\times K)$, $l,k\in \N$, for a set $K$ is the space of real valued functions $\varphi$ of class $C^l$ w.r.t. the first variable and measurable w.r.t. the second.  For $(z,k)\in Z\times K$, we denote by $D^{l}\varphi(z,k)$ the (Fr\'echet) differential of $\varphi$, of order $l$ w.r.t. $z$, computed at $(z,k)$.
A subscript $b$ can be added - to obtain $C_b^{l,0}(Z\times K)$ - in order to specify that the functions and their succesive differentials are uniformly bounded. 

$\disp D(I)$ is the space  of right-continuous, left-limited (or \textit{c\`adl\`ag}\footnote{From French \textit{continu \`a droite et admettant une limite \`a gauche}.}) real valued functions defined on $I=[0,1]$, it is endowed with the Skorohod topology;
 $\disp C(I)$ (resp. $C^k(I)$) is the space of periodic continuous (resp. $C^k$) real valued functions defined on $I=[0,1]$;
 $\disp C\big([0,T],Z\big)$ (resp. $\disp C\big(\R^+,Z\big)$) is the space of continuous processes defined on $[0,T]$ (resp. $\R^+$) with values in $Z$; 
$\disp D\big([0,T],Z\big)$ (resp. $\disp D\big(\R^+,Z\big)$ is the space cadlag processes defined on $[0,T]$ (resp. $\R^+$) with values in $Z$. It is endowed with the Skorohod topology. 

\subsection{Infinite dimensional Piecewise Deterministic Processes}\label{infinite_dimensional_PDMPs}

Piecewise deterministic Markov processes (PDMPs) form a class of processes that has been formalized in the finite dimensional case in \cite{Davis1984,Davis1993} among others. We briefly define infinite dimensional (continuous) PDMPs and give some useful results about them. We refer to \cite{Riedler2011} or \cite{Riedler2011bis} for an in depth presentation.

Consider a Banach space $B$ endowed with a norm $\Vert\cdot\Vert_B$, a countable set of isolated states $K$ equipped with the discrete topology, and set $E:=B\times K$. Then, let $\disp \big(\Omega,\fF,(\fF_t)_{t\geq0},\P\big)$ denote a filtered probability space satisfying the usual conditions.\par 

In this article, a standard "continuous" PDMP with values in $E$ is a non exploding c\`adl\`ag stochastic process $\disp \big\{ u(t)=(u_C(t),u_D(t)), t\geq0\big\}$, determined by its four characteristics:\vspace{0.2cm}\par 

\textbf{(1)-(2)}\textit{Linear and nonlinear operators:} for every $\nu\in K$, there is given an abstract evolution equation
\begin{equation}\label{abstract_evolution_equation_of_a_PDMP}
\frac{d}{dt} u_C(t)=L_\nu u_C(t)+F_\nu(u_C(t)),
\end{equation}
where $F_\nu: B\rightarrow B$ is a (possibly nonlinear) operator and $L_\nu: B\rightarrow B$ is an unbounded linear operator, which is $m-$dissipative with dense domain. Let $\{S_\nu(t):=e^{L_\nu t},t\geq0\}$ be the semigroup generated by $L_\nu$. We assume that there exists a unique global mild flow $\phi_\nu(\hspace{0.1cm}\cdot\hspace{0.1cm},\alpha)$ to $(\ref{abstract_evolution_equation_of_a_PDMP})$, satisfying $\disp \phi_\nu(\hspace{0.1cm}\cdot\hspace{0.1cm},\alpha)\in C\big(\R_+,B\big)$ and \[ \phi_\nu(t,\alpha)=S_\nu(t)\alpha+\int_0^tS_\nu(t-s)F_\nu\big(\phi_\nu(s,\alpha)\big)ds, \] for periodic boundary conditions, every initial value $\alpha\in B$, and all $\nu\in K$. We often use the notations $L_\nu(\alpha)=L(\alpha,\nu)$, $F_\nu(\alpha)=F(\alpha,\nu)$ and $\disp \phi_\nu(\hspace{0.1cm}\cdot\hspace{0.1cm},\alpha)=\phi(\hspace{0.1cm}\cdot\hspace{0.1cm},\alpha,\nu)$.\vspace{0.1cm}\par 

\textbf{(3)}\textit{A jump rate:} $\Lambda: E\rightarrow\R_+$, which is measurable, and such that for all $(\alpha,\nu)\in E$, the function $t\mapsto\Lambda\big(\phi_\nu(t,\alpha),\nu\big)$ is integrable over every finite time interval, but divergent over $\R_+$. In other words, 
for all $(\alpha,\nu)\in E$ and $T>0$, \[ \int_0^{T}\Lambda\big(\phi_\nu(t,\alpha),\nu\big)dt<\infty\hspace{0.5cm}\text{while}\hspace{0.5cm}\int_0^\infty\Lambda\big(\phi_\nu(t,\alpha),\nu\big)dt=\infty. \] Also, we often use the notation $\Lambda_\nu\big(\phi_\nu(t,\alpha)\big)=\Lambda\big(\phi_\nu(t,\alpha),\nu\big)$. \vspace{0.1cm}\par 

\textbf{(4)}\textit{A transition measure:} $\displaystyle Q : E \longrightarrow \mathcal{P}\left(K\right)$ which is measurable, such that for every fixed $A\subset K$, the function $\disp (\alpha,\nu)\mapsto Q_{\nu}(A;\alpha)$ is measurable, and satisfies $\disp Q_\nu\big(\{\nu\};\alpha\big)=0$ for all $(\alpha,\nu)\in E$. Again, we often denote $\disp Q_{\nu}(A;\alpha)=Q(A;\alpha,\nu)$.

The second component $u_D$ of the PDMP takes (discrete) values in $K$, and is called the discrete component. It has right continuous piecewise constant sample paths and is often denoted the piecewise constant or jump component. The first component $u_C$ takes (continuous) values in $B$, and is called the continuous component. It has continuous sample paths and justifies the name "continuous" PDMP. The mechanism which governs the evolution of the process is as follows. While $u_D$ is constant with $u_D=\nu$, $u_C$ evolves according to the flow determined by the operators $L_\nu$ and $F_\nu$ through (\ref{abstract_evolution_equation_of_a_PDMP}). Then, a jump occurs at a random time with the jump rate $\Lambda$, and the target state after that jump is determined by the transition measure $Q$.\\

\noindent \textbf{\textit{Construction.}} A \textit{c\`adl\`ag} sample path $\disp \big\{v(t),t\geq0\big\}$ starting at 
\[ v(0)=\big(v_C(0),v_D(0)\big)=(\alpha,\nu)=v_0\in E \] can be constructed for such a process as follows. \vspace{0.1cm}\par 

$\bullet$ For all $t\in [0,T_1)$, \[ v(t)=\big(\phi_\nu(t,\alpha),\nu\big), \] where $T_1$ is the first jump time of $v_D$. Denote by $\tau_1$ the waiting time of the first transition. Then $\tau_1=T_1$, and, conditionally to the starting point, $\tau_1$ has an exponential distribution. Its survivor function $H_\nu:\R_+\times B\rightarrow \R_+$, is defined by 
\begin{equation}\label{survivor_function_of_the_first_transition_time}
\disp H_\nu(t,\alpha):=\P_{v_0}\big\{\tau_1>t\big\}=\text{exp}\left(-\int_0^t\Lambda\big(\phi_\nu(s,\alpha),\nu\big)ds\right).
\end{equation}
It is often convenient to use the notation $H(t,\alpha,\nu)=H_\nu(t,\alpha)$.

$\bullet$ At time $T_1^-$, \[ v(T_1^-)=\big(\phi_\nu(T_1,\alpha),\nu\big). \]

$\bullet$ A transition occurs at time $T_1$. The target state $v(T_1)=\big(v_C(T_1),v_D(T_1)\big)$ satisfies 
\[ v_C(T_1)=\phi_\nu(T_1,\alpha)=v_C(T_1^-)=\alpha_1, \] 
and $v_D(T_1)$ has the distribution \[ \P_{v_0}\big\{v_D(T_1)\in A\big|T_1=t\big\}=Q_\nu\big(A;\phi_\nu(t,\alpha)\big), \] 
for every measurable subset $A\subset K$, and $t>0$.\par 

$\bullet$ After that first transition, the jump component $v_D$ remains constant and equal to the chosen target state $v_D(T_1)=\nu_1$, until its next jump. The continuous component $v_C$ then evolves according to the "updated" abstract evolution equation, starting from $v_C(T_1)=\alpha_1$ with $v_D=\nu_1$. The procedure is repeated independently starting from $(\alpha_1,\nu_1)$, and the process is recursively constructed. One obtains a sequence $\tau_1,\tau_2,\cdots$ of independent transition waiting times and a sequence of jump times $T_1,T_2,\cdots$, with $T_i=\tau_1+\cdots+\tau_i$ for all $i\geq1$. \vspace{0.1cm}\par 

The constructed process is a "continuous" PDMP. The number of jumps that occur between times $0$ and $t$ is 
\[ N_t=\sum_{i\geq1}\1_{(T_i\leq t)}. \] 
The non explosion of the process is usually refered to, as the regularity of the PDMP. It means an almost sure finite number of jumps until finite times, and is characterized by $T_i\rightarrow\infty$ as $i\rightarrow\infty$. However, that condition is not easy to check in practice. Still, it is satisfied in particular when the expected number of jumps is finite on any finite time interval (see \cite{Davis1993}, p.60). Therefore, in order to make sure our PDMP is regular, we assume, as part of the standard conditions, that
\begin{assumption}\label{non_explosion_of_the_PDMP}
$\E_{(\alpha,\nu)}\big[N_t\big]<\infty,\hspace{0.5cm}\forall (\alpha,\nu)\in E,\hspace{0.1cm} t\geq0$.
\end{assumption}

Next, a PDMP characterized by a quadruple $(L,F,\Lambda,Q)$ can be constructed as above in some probability space, to be a strong Markov \textit{c\`adl\`ag} process (\cite{Riedler2011bis}, Theorem 4). It is then called a PDMP, meaning piecewise deterministic Markov process. Only standard PDMPs are considered in the sequel.\\

\noindent \textbf{\textit{The full generator and the martingale probem.}} Consider a PDMP $v=(v_C,v_D)$ with the characteristics $(L,F,\Lambda,Q)$ and let $(P_t)\equiv (P_t)_{t\geq0}$ be the corresponding semigroup on $E=B\times K$. 
For all $t\geq0$, $P_t:\fB_b(E)\rightarrow \fB_b(E)$, and \[ P_t\varphi(\alpha,\nu)=\E_{(\alpha,\nu)}\big[\varphi(u(t))\big] \] for all $\varphi\in \fB_b(E)$, $(\alpha,\nu)\in E$. 
We have denoted by $\E_{(\alpha,\nu)}$ the conditional expectation given $v(0)=(\alpha,\nu)$. 
Clearly, $(P_t)$ is a semigroup of contraction and is measurable in the sense that the map $t\mapsto P_t\varphi$ is measurable for all $\varphi\in\fB_b(E)$. 
Associated with $(P_t)$ is its full generator 
\[ \hat{\cA}=\left\{(\varphi,\psi) \in\fB_b(E)\times\fB_b(E):P_t\varphi-\varphi=\int_0^tP_s\psi ds\right\}. \] 
This operator is in general multi-valued, i.e. given any $\varphi$, there needs not be a unique $\psi$ such that $(\varphi,\psi)\in\hat{\cA}$. Hence, for a given semigroup, the associated full generator can not be viewed as an operator on $\fB_b(E)$ in general. Find more details in \cite{Kurtz1986}, Part 1, Section 5. However, we see below that the situation is simpler for PDMPs.

Now, denote by $P_{v_0}$ the law of the PDMP $v$ when it starts from $v_0\in E$. 
It is well known that $P_{v_0}$ is a solution of the martingale problem associated with $\hat{\cA}$ in the sense:
 \[ \varphi(w(t))-\varphi(v_0)-\int_0^t\psi(w(s))ds \] defines a $P_{v_0}$-martingale, for all $(\varphi,\psi)\in \hat{\cA}$. 
We have denoted by $\disp \big\{w(t),t\geq0\big\}$ the canonical process on the probability space $\disp \big(D(\R_+,E),\cD_E,P_{v_0}\big)$. 
An equivalent formulation is: 
\begin{equation}\label{martingale_characterising_the_full_generator}
\disp N_\varphi(t):=\varphi(v(t))-\varphi(v_0)-\int_0^t\psi(v(s))ds 
\end{equation} 
defines a $\P_{v_0}$-martingale, for all $(\varphi,\psi)\in \hat{\cA}$. One is easily convinced that the full generator $\hat{\cA}$ is exactly the subset of all the couples $(\varphi,\psi)\in\fB_b(E)\times\fB_b(E)$ for which the martingale problem above is satisfied. \\



\noindent \textbf{\textit{The extended generator.}} We notice that the martingale problem is an essential tool for characterizing the law of a PDMP. However, it is directly related to the full generator of the process, which in turn is very often not easy to determine explicitely. 
With that idea in mind, the extended generator associated with the PDMP is considered, which is the operator $\bar{\cA}$, whose domain is 
\[ \cD(\bar{\cA}):=\big\{\varphi\in\fB(E), \exists \psi\in\fB(E) : N_\varphi(t) \text{ defines a } \P_{v_0}-\text{local martingale}\big\}, \]
where again
\begin{equation}\label{local_martingale_characterising_the_extended_generator}
\disp N_\varphi(t):=\varphi(v(t))-\varphi(v_0)-\int_0^t\psi(v(s))ds.
\end{equation} 
Clearly, it is an extension of the full generator.  This justifies the name "extended" generator. 

It should be emphasized, see \cite{Davis1993} pp 32-33, that the extended generator is a single-valued operator on $\fB(E)$, up to sets of \textit{zero potential}, these are sets $A\subset\fB(E)$ such that 
\[ \int_0^\infty\1_A(v(s))ds=0\hspace{1cm}\P_{v_0}-a.s.\hspace{0.5cm}\text{for every }\hspace{0.2cm}v_0\in E. \] 
The process "spends no time" in $A$, regardless of the starting point. Thus, we can set: $\bar{\cA}\varphi=\psi$.

Therefore, the full generator is in particular also single-valued up to sets of \textit{zero potential}. 
Thus, 
it may be considered as an operator $\hat{A}$ on $\fB(E)$, with the domain 
\[ \cD(\hat{\cA})=\left\{\varphi\in\fB(E):P_t\varphi-\varphi=\int_0^tP_s\hat{\cA}\varphi ds\right\}. \]
Furthermore, if the test function $\varphi$ is such that $N_\varphi(t)$ given by \eqref{local_martingale_characterising_the_extended_generator} is bounded, then $N_\varphi(t)$ defines a martingale and hence, the restriction of $\bar{\cA}$ to such test functions coincides with $\hat{\cA}$.\par 

The domain of the extended generator is characterized in \cite{Davis1993}, Theorem 26.14, for finite dimensional PDMPs. The infinite 
dimensional case is considered  in \cite{Riedler2011bis}, Theorem 4. 
Following their arguments, we know that, for well chosen test functions $\varphi$,  the generator may be identified with the operator  
\begin{equation}\label{extended_generator_of_a_PDMP}
\begin{array}{l}
\disp \cA\varphi(\alpha,\nu)=\big\langle D^{1,0}\varphi(\alpha,\nu), L(\alpha,\nu)+F(\alpha,\nu)\big\rangle\\
\disp\hspace{3.5cm} +\Lambda(\alpha,\nu)\int_K\big[\varphi(\alpha,\xi)-\varphi(\alpha,\nu)\big]Q(d\xi,\alpha,\nu),
\end{array}
\end{equation}
where $\langle\cdot,\cdot\rangle$ is the duality paring between $B$ and its topological dual $B'$, and $D^{1,0}\varphi(\alpha,\nu)$ the 
Fr\'echet differential of $\varphi$ w.r.t. its first variable. 

In the case when there is an underlying Hilbert space $H$ such that $B\subset H$, one identifies $D^{1,0}\varphi(\alpha,\nu)$ with the corresponding gradient, and the duality paring with the inner product $\langle\cdot,\cdot\rangle$ of $H$. 
The main difference with the finite dimensional case is that the expression above contains an unbounded operator and in general $L(\alpha,\nu)$ is not in $B$. We need to restrict the test functions so that the duality pairing $\big\langle D^{1,0}\varphi(\alpha,\nu), L(\alpha,\nu)\big\rangle$ is meaningful. Another possibility 
is to restrict the values of $\alpha$ to a smaller space than $B$.

Below, we identify a subset $\cE$ of the domain of $\cA$ such that for $\varphi\in\cE$ \cref{extended_generator_of_a_PDMP} holds.

\subsection{Multiscale stochastic spatial regulatory networks with a slow dynamic independent of the space discretization}
\label{s2.3}

Molecules of two species $C$ and $D$ are submitted to reactions and diffusions in the unit interval, the spatial domain $I$. Following \cite{Arnold1980}, we devide $I$ into $N$ smaller intervals, called \textit{sites}, of equal length $\displaystyle N^{-1}$: $\displaystyle I_j=\left((j-1)N^{-1},jN^{-1}\right]$, for $j=1,\cdots,N$. Molecules are produced (birth) or removed (death) on sites at rates which depend on the local current number of particles. Moreover, molecules of $C$ can diffuse between sites by simple random walks (one at once), with jump rates proportional to $N^2$ and linearly depending on the current local state. In this framework, an event can be either an onsite chemical reaction or a diffusion. We often say "reaction" for "event", and use the super/subscript $C$ (resp. $D$) in reference to $C$ (resp. $D$).

The species $C$ has a large population size scale while that of $D$ is small. Molecules of $D$ do not diffuse. As in \cite{Arnaud2012}, we divide the set $\mathfrak{R}$ of possible onsite reactions in three disjoint subsets: 
\[ \mathfrak{R}=\mathfrak{R}_C\cup\mathfrak{R}_{DC}\cup\mathfrak{R}_D. \] 
Reactions in $\fR_C$ (resp. $\fR_D$) involve only reactants and products of type $C$ (resp. $D$), whereas, reactions in $\mathfrak{R}_{DC}$ involve both types of reactants and/or products. Also, \\
 
$\bullet$ $\displaystyle\hspace{0.3cm} X_j^{N,C}$\hspace{0.1cm} $\displaystyle \left(\hspace{0cm} \text{resp.} \hspace{0.1cm}  X_j^{N,D}\hspace{0.1cm}\right)$ is the number of molecules of $C$ (resp. $D$) on the site $j$; \vspace{0.1cm}\par 
 
$\bullet$ $\displaystyle\hspace{0.3cm} X_C^N:=\left(X_j^{N,C}\right)_{1\leq j\leq N}$\hspace{0cm},\hspace{0.2cm} $\displaystyle X_D^N:=\left(X_j^{N,D}\right)_{1\leq j\leq N}$\hspace{0.2cm} and \hspace{0.2cm}$\displaystyle X^N:=\left(X_C^N,X_D^N\right)\in\N^{2N}$.\\ 

\noindent  For simplicity, we consider periodic boundary conditions: $\displaystyle X_{j+N}^N=X_{j}^N$, $\forall j$. Hence, the molecular composition of the system $\displaystyle X^N:=\left\{X^N(t),t\geq0\right\}$ is a $\N^{2N}$-valued Markov process with the transitions: 
\[
\left\{
\begin{array}{l}
\displaystyle \left(X_C^N,X_D^N\right)\longrightarrow \left(X_C^N+\gamma_{j,r}^{C}e_j,X_D^N\right)\hspace{0.2cm} \text{at rate}\hspace{0.1cm} \lambda_r\left(X_j^{N,C}\right),\hspace{0.1cm} \text{for}\hspace{0.1cm} r\in\mathfrak{R}_C,\vspace{0.1cm}\\

\displaystyle \left(X_C^N,X_D^N\right)\longrightarrow \left(X_C^N+\gamma_{j,r}^{C}e_j,X_D^N+\gamma_{j,r}^{D}e_j\right)\hspace{0.2cm} \text{at rate} \hspace{0.1cm} \lambda_r\left(X_j^{N,C},X_j^{N,D}\right),\hspace{0.1cm} \text{for}\hspace{0.1cm} r\in\mathfrak{R}_{DC},\vspace{0.1cm}\\

\displaystyle \left(X_C^N,X_D^N\right)\longrightarrow \left(X_C^N,X_D^N+\gamma_{j,r}^{D}e_j\right)\hspace{0.2cm} \text{at rate}\hspace{0.1cm} \lambda_r\left(X_j^{N,D}\right),\hspace{0.1cm} \text{for}\hspace{0.1cm} r\in\mathfrak{R}_D,\vspace{0.1cm}\\

\displaystyle \left(X_C^N,X_D^N\right)\longrightarrow \left(X_C^N+e_{j-1}-e_j,X_D^N\right)\hspace{0.2cm} \text{at rate}\hspace{0.1cm} N^2X_{j}^{N,C}, \hspace{0.1cm} \text{for a diffusion}\hspace{0.1cm}j\rightarrow j+1,\vspace{0.1cm}\\

\displaystyle \left(X_C^N,X_D^N\right)\longrightarrow \left(X_C^N+e_{j+1}-e_j,X_D^N\right)\hspace{0.1cm} \text{at rate}\hspace{0.1cm} N^2X_{j}^{N,C}, \hspace{0.1cm} \text{for a diffusion}\hspace{0.1cm}j\rightarrow j-1,
\end{array}
\right.
\]
\noindent where $\left\lbrace e_j, j=1,\cdots,N\right\rbrace$ is the canonical basis of $\mathds{R}^N$, and $\gamma_{j,r}^{C},\gamma_{j,r}^{D}\in\mathds{Z}$, $\forall r\in\mathfrak{R}$, $\forall 1\leq j\leq N$.\\

\noindent \textbf{\textit{Scaling and Density dependence.}} On every site, the initial average number of molecules for $C$ is of order $\mu$ with $\mu$ large. Namely, if $\disp \mathfrak{M}=\mathfrak{M}_C+\mathfrak{M}_D$ is the total initial number of molecules, then $\displaystyle \mathfrak{M}_C\approx N\times\mu$.  The size of $\disp \mathfrak{M}_D$ is precised below.

Then, we make the

\begin{assumption}\label{density_dependence_2}\hspace{1cm}\par 
\textbf{(i)} Density dependence holds for the rates of reactions in $\fR_C$. That is, for all $r\in\fR_C$, there exists $\tilde{\lambda}_r$ satisfying \[ \lambda_r\left(X_j^{N,C}\right)=\mu\tilde{\lambda}_r\left(\frac{X_j^{N,C}}{\mu}\right)\hspace{0.4cm} \text{for all}\hspace{0.3cm} j=1,\cdots,N. \]\par 
\textbf{(ii)} Reactions $r\in\fR$ are spatially homogeneous: $\gamma_{j,r}^C=\gamma_r^C$, $\gamma_{j,r}^D=\gamma_r^D$.
\par 
\textbf{(iii)} The molecules of $C$ diffuse, while those of $D$ do not. 
\end{assumption}

It follows that reactions in $\fR_C$ are fast, while reactions in $\fR_D$ are slow. Below, we distinguish two types of mixed reactions in $\fR_{DC}$, some will be fast. 

For notational convenience, we omit the tilde for reaction rates obtained by density dependence, as described in \cref{density_dependence_2} (i). 
Moreover, we consider the parameter $\mu$ as a function of $N$ which goes to infinity with $N$, and omit to mention the dependance on $\mu$.
The assumption of spatial homogeneity is not essential. It allows simpler notations.

Then, we rescale $X^N$ and define 
 $$  U_j^{N,C}=\frac{X_j^{N,C}}{\mu}, \; U_j^{N,D}=X_j^{N,D};\; U_C^N=\left(U_j^{N,C}\right)_{1\leq j\leq N},\; U_D^N=\left(U_j^{N,D}\right)_{1\leq j\leq N},$$
 and
$$ U^N= \left(U_C^{N},U_D^{N}\right)\in\mathds{R}^{N}\times\N^{N}.$$ 
Since $U_C^N$ (resp. $U_D^N$) has continuous (resp. discrete) values, $C$ (resp. $D$) is said to be continuous (resp. discrete). The generator of the new scaled process $U^N$ is given by\\

\noindent $\displaystyle \mathcal{A}_0^N\varphi(U):=\left.\frac{d}{dt}\mathds{E}_U\left[\varphi\left(U^N(t)\right)\right]\right|_{t=0}\hspace{0.2cm}:=\hspace{0.2cm}\lim_{t\rightarrow0}\frac{1}{t}\mathds{E}_U
\left[\varphi\left(U^N(t)\right)-\varphi(U)\right]$\par 

\noindent $\displaystyle \hspace{0.5cm} =\sum_{j=1}^N\left\{\sum_{r\in\mathfrak{R}_C}\left[\varphi
\left(U_C+\frac{\gamma_{r}^{C}}{\mu}e_j,U_D\right)-\varphi(U_C,U_D)\right]\mu\lambda_r\left(U_j^C\right)\right.$\par 

$\disp \hspace{1.7cm}+\sum_{r\in\mathfrak{R}_{DC}}\left[\varphi
\left(U_C+\frac{\gamma_{r}^{C}}{\mu}e_j,U_D+\gamma_r^{D}e_j\right)-\varphi(U_C,U_D)\right]\lambda_r\left(U_j^C,U_j^D\right)$\par 

$\disp \hspace{1.7cm}+\left.\sum_{r\in\mathfrak{R}_{D}}\left[\varphi
\left(U_C,U_D+\gamma_r^{D}e_j\right)-\varphi(U_C,U_D)\right]\lambda_r\left(U_j^D\right)\right\}$\par 

\noindent $\displaystyle \hspace{0.3cm}+\sum_{j=1}^N\left\{\left[\varphi\left(U_C+\frac{e_{j-1}-e_j}{\mu},U_D\right)+\varphi\left(U_C+\frac{e_{j+1}-e_j}{\mu},U_D\right)-2\varphi(U_C,U_D)\right] \mu N^2U_j^C\right\}$,\\

\noindent on the domain $\displaystyle \fB_b\left(\mathds{R}^{2N}\right)$ (see \cite{Kurtz1971} or \cite{Kurtz1986} pp162-164). 

In order to achieve a pointwise modeling over the whole spatial domain, we introduce the step function 
\begin{equation}\label{time_space_function2}
u^N(t,x)=\sum_{j=1}^NU_j^N(t)\mathds{1}_j(x), \hspace{0.3cm} t\geq 0,\hspace{0.2cm} x\in I_j,
\end{equation} 
\noindent where $\displaystyle \mathds{1}_j(\cdot):=\mathds{1}_{I_j}(\cdot)$ is the indicator function of the $j-$th site $\displaystyle I_j$. Note that for all $t\geq0$, the function $u^N(t)=u^N(t,\cdot)$ can be identified with the vector $U^N(t)$ of $\R^{2N}$. It is a $1-$periodic function, since $U_{j+N}^N=U_j^N$ for all $j$. Now, let $\displaystyle \mathds{H}^N$ denote the subspace of $L^2(I)$ which consists in real-valued step functions defined on $I$, and constant on every site $I_j$, $1\leq j\leq N$. We extend functions in $\displaystyle \mathds{H}^N$ to be periodic. Moreover, consider the canonical projection 
\begin{equation}\label{Canonical_Projection2}
\begin{tabular}{lll}
$\displaystyle P_N$ & $:$ & $L^2(I) \longrightarrow \mathds{H}^N$\vspace{0cm}\\
 & & $\displaystyle u\longmapsto P_Nu=\sum_{j=1}^Nu_j\mathds{1}_j,\hspace{0.3cm}\text{with}\hspace{0.3cm} u_j:=N\int_{I_j}u(x)dx$.
\end{tabular}
\end{equation}
\noindent Then,\par 

$\disp \hspace{1.3cm} u_C^N(t,x)=\sum_{j=1}^Nu_j^{N,C}(t)\mathds{1}_j(x)\hspace{0.5cm}\text{and}\hspace{0.4cm}u_D^N(t,x)=\sum_{j=1}^Nu_j^{N,D}(t)\mathds{1}_j(x)$,\par 

\noindent where\par 

$\disp\hspace{2cm}u_j^{N,C}(t):=P_Nu_C^N(t)=N\int_{I_j}u_C^N(t,x)dx = U_j^{N,C}(t)$,\\

\noindent and a similar relation holds for $D$. The process, $\displaystyle u^N:=\left\{u^N(t),t\geq 0\right\}$ is a $\displaystyle \mathds{H}^N\times\mathds{H}^N$-valued \textit{c\`adl\`ag} 
Markov process, with the transitions:
\begin{equation}\label{transitions_for_the_spatial_model2}
\left\{
\begin{tabular}{l}
$\displaystyle \left(u_C^N,u_D^N\right)\longrightarrow\left(u_C^N+\frac{\gamma_r^C}{\mu}\mathds{1}_j,u_D^N\right)$,\hspace{0cm} at rate $\displaystyle \mu\lambda_r\left(u_j^{N,C}\right)$, \hspace{0cm} for $r\in\mathfrak{R}_C$,\vspace{0.1cm}\\

$\displaystyle \left(u_C^N,u_D^N\right)\longrightarrow\left(u_C^N+\frac{\gamma_{r}^{C}}{\mu}\1_j,u_D^N+\gamma_r^{D}\1_j\right)$,\hspace{0cm} at  rate $\displaystyle \lambda_r\left(u_j^{N,C},u_j^{N,D}\right)$,\hspace{0cm} $r\in\mathfrak{R}_{DC}$,\vspace{0.1cm}\\

$\displaystyle \left(u_C^N,u_D^N\right)\longrightarrow\left(u_C^N,u_D^N+\gamma_r^{D}\1_j\right)$,\hspace{0cm} at rate $\displaystyle \lambda_r\left(u_j^{N,D}\right)$,\hspace{0cm} for $r\in\mathfrak{R}_{D}$,\vspace{0.1cm}\\

$\displaystyle \left(u_C^N,u_D^N\right)\longrightarrow \left(u_C^N+\frac{\mathds{1}_{j-1}-\mathds{1}_{j}}{\mu},u_D^N\right)$,\hspace{0cm} at rate $\displaystyle \mu N^2u_j^{N,C}$,\vspace{0.1cm}\\

$\displaystyle \left(u_C^N,u_D^N\right)\longrightarrow\left(u_C^N+\frac{\mathds{1}_{j+1}-\mathds{1}_{j}}{\mu},u_D^N\right)$,\hspace{0cm} at rate $\displaystyle  \mu N^2u_j^{N,C}$.
\end{tabular}
\right.
\end{equation}
\noindent 
Such a Markov process does exist and is unique (see \cite{Kotelenez1988bis}, which is based on \cite{Kurtz1986}) until a possible blow-up time. 
%
In addition, under natural assumptions on the reaction rates, we have $u^N(t)\geq0$ for all $t$, as soon as $u^N(0)\geq0$. 
\vspace{0.2cm}\par 

We now specify the description of the mixed reactions $r\in\mathfrak{R}_{DC}$.

\begin{assumption}\label{partitioning_the_set_of_mixed_reactions2} In some $\displaystyle S_1\subset\mathfrak{R}_{DC}$, reactions are spatially homogeneous, fast and do not affect the discrete species:\par 
$\bullet$ $\gamma_{r}^D=0$ and the rate is $\lambda_r\left(u_j^{N,C},u^{N,D}_j\right)=\mu\tilde\lambda_r\left(u_j^{N,C},u^{N,D}_j\right)$, for $r\in S_1$ (\textit{fast mixed reactions}).

\end{assumption}
Again, below we omit the tildes. 

\medskip
This model leads to mathematical difficulties as explained in \cite{DebusscheNankep2017}. 
Unless $u_{j}^{N,D}$ is zero for almost all $j$, the discrete species are in fact also abundant since the total number is of order $N$, hence they also converge to a continuous model. The 
limit model is expected to be a PDE for $u_C$ coupled to an ODE for $u_D$. As explained in \cite{DebusscheNankep2017}, the difficulty is that the convergence seems to hold in bad topologies
and this limit is probably difficult to justify for nonlinear reaction rates. Thus, in \cite{DebusscheNankep2017}, a spatial correlation in the reaction involving the discrete species has been introduced and the limit model is 
described by a PDE coupled to a nonlocal differential equation. 

\medskip
In this paper, we introduce another situation where the discrete species remain discrete at the limit. We introduce an assumption of spatial multiscaling: the spatial domain $I$ is split into a finite number 
of macrosites $J_1,\dots,J_k$, where $k$ is a fixed finite number, and $u_D$ is constant in each $J_\ell$. 
This describes a situation where the discrete species are of much larger size than the continuous species. For instance the discrete species are constant in alls cells, while the continuos ones are varying inside each cell. Similarly, the multi-scaling could be between groups of cells and cells

\medskip In order to avoid complicated notations, we assume for simplicity that all $J_\ell$ has the same length and consider $N$ as multiples of $k$. Thus each $J_\ell$ is the union of microsites:
$$
J_\ell=\bigcup_{j=N_{\ell-1}+1}^{N_\ell} I_j,\; \ell= 1, \dots,k,
$$
with $N_\ell=\frac{\ell N}k$. 
On each $J_\ell$, $u^N_D$ is constant and takes the value $u^{N,D}_\ell$:
$$
u^N_D(t,x)=\sum_{\ell=1}^k u^{N,D}_\ell(t) \mathds{1}^\ell(x),
$$
where $\mathds{1}^\ell=\mathds{1}_{J_\ell}$ is the indicator function of the $\ell-$th macrosite $J_\ell$. 

It is then natural to consider that when there is a slow reaction between discrete and continuous species, {\it i.e.}  in $\fR_{DC}\backslash S_1$,  on $J_\ell$, it may affect $u_C$ on all $J_\ell$:
$$
\displaystyle \left(u_C^N,u_D^N\right)\longrightarrow\left(u_C^N+\frac{\gamma_{r}^{C}}{\mu}\sum_{j=N_{\ell-1}+1}^{N_\ell}b_{j,N}^r\1_j ,u_D^N+\gamma_r^{D}\1^\ell\right),
$$ 
for some non negative $b_{j,N}^r$. Below we use the notation:
$$
\gamma_{r,\ell}^{C,N} =\gamma_{r}^{C}\sum_{j=N_{\ell-1}+1}^{N_\ell}b_{j,N}^r\1_j.
$$  
We choose $b^r_{j,N}= N \int_{I_j} b^r(x)dx$ for some positive function $b^r\in C([0,1])$. 

We also assume that the rate of such reaction depends on $u_\ell^{N,D}$, and on 
$u_C^N$ (considered all over $J_\ell$) through a local average. In other words, the rate of such reaction is of the form $\lambda_r\left(\sum_{j=N_{\ell-1}+1}^{N_\ell}a_{j,N}^r u_j^{N,C},u^{N,D}_\ell\right)$ where $a_{j,N}^r$ 
are non negative numbers summing to $1$: $\sum_{j=N_{\ell-1}+1}^{N_\ell}a_{j,N}^r=1$. We choose $a^r_{j,N}= \int_{I_j} a^r(x)dx$ for some positive function $a^r\in C([0,1])$ such that $\int_0^1 a^r(x)dx=1$.

\medskip

In this framework ,  \eqref{transitions_for_the_spatial_model2} becomes
\begin{equation}\label{transitions_for_the_spatial_model3}
\left\{
\begin{tabular}{l}
$\displaystyle \left(u_C^N,u_D^N\right)\longrightarrow\left(u_C^N+\frac{\gamma_r^C}{\mu}\mathds{1}_j,u_D^N\right)$,\hspace{0cm} at rate $\displaystyle \mu\lambda_r\left(u_j^{N,C}\right)$, \hspace{0cm} for $r\in\mathfrak{R}_C$,\vspace{0.1cm}\\

$\displaystyle \left(u_C^N,u_D^N\right)\longrightarrow\left(u_C^N+\frac{\gamma_{r}^{C}}{\mu}\1_j,u_D^N\right)$,\hspace{0cm} at  rate $\displaystyle \mu\lambda_r\left(u_j^{N,C},u_{\ell_j}^{N,D}\right)$,\hspace{0cm} for $r\in S_1$,\vspace{0.1cm}\\

$\displaystyle \left(u_C^N,u_D^N\right)\longrightarrow\left(u_C^N+\frac{\gamma_{r,\ell}^{C,N}}{\mu} ,u_D^N+\gamma_r^{D}\1^\ell\right)$,\hspace{0cm} at  rate $\displaystyle \lambda_r\left(\sum_{j=N_{\ell-1}+1}^{N_\ell}a_{j,N}^ru_j^{N,C},u^{N,D}_\ell\right)$,\hspace{0cm}\\
\hfill for $r\in\mathfrak{R}_{DC}\backslash S_1$,\vspace{0.1cm}\\

$\displaystyle \left(u_C^N,u_D^N\right)\longrightarrow\left(u_C^N,u_D^N+\gamma_r^{D}\1^\ell\right)$,\hspace{0cm} at rate $\displaystyle \lambda_r\left(u_\ell^{N,D}\right)$,\hspace{0cm} for $r\in\mathfrak{R}_{D}$,\vspace{0.1cm}\\

$\displaystyle \left(u_C^N,u_D^N\right)\longrightarrow \left(u_C^N+\frac{\mathds{1}_{j-1}-\mathds{1}_{j}}{\mu},u_D^N\right)$,\hspace{0cm} at rate $\displaystyle \mu N^2u_j^{N,C}$,\vspace{0.1cm}\\

$\displaystyle \left(u_C^N,u_D^N\right)\longrightarrow\left(u_C^N+\frac{\mathds{1}_{j+1}-\mathds{1}_{j}}{\mu},u_D^N\right)$,\hspace{0cm} at rate $\displaystyle \mu N^2u_j^{N,C}$,
\end{tabular}
\right.
\end{equation}
where here and below $\ell_j$ is the upper integer part of $\frac{jk}{N}$. In other words: $N_{\ell_j-1}+1\le j\le N_{\ell_j}$. 

\medskip
Note that this system has a problem since it does not ensure positivity of $u_C^N$. In general, this is ensured by the fact that the rate vanishes when the concentration vanishes. But this is not possible for the type of rate considered for the mixed reactions in $\mathfrak{R}_{DC}\backslash S_1$. A more realistic rate for such reactions would be:
$$
\displaystyle \lambda_r\left(\sum_{j=N_{\ell-1}+1}^{N_\ell}a_{j,N}^ru_j^{N,C},u^{N,D}_\ell\right)\prod_{j=N_{\ell-1}+1}^{N_\ell}\1_{u_j^{N,C}+\frac{\gamma^r_C b_{j,N}^r}{\mu}\ge 0}
$$
When $\mu\to \infty$, this further factor converges to $\prod_{j=N_{\ell-1}+1}^{N_\ell}\1_{u_j^{N,C}\ge 0}$ or $\prod_{j=N_{\ell-1}+1}^{N_\ell}\1_{u_j^{N,C}> 0}$ depending on the sign of $\gamma^r_C$. This extra factor has no effect at the limit since $u_C$ is positive for $\mu=\infty$. Indeed, positivity is ensured by the other reactions acting on the continuous component.

Moreover, due to the lack of smoothness of the indicator function, we should replace $\1_{u\ge 0}$ by a smooth function approximating it. This creates further notational complexity but no mathematical problems (see \cite{DebusscheNankep2017}) and we do 
not consider such extra factor.
\medskip

We follow the spatial distribution of the discrete component at the level of the macrosites $J_\ell$.
We may identify the space of such functions with $\R^k$ thanks to the formula: $u_D(t,x)=\sum_{\ell=1}^k u_{\ell}^D(t)\1^\ell(x)$.
%
%
Thus, with this new framework, the phase space is $\mathds{H}^{N}\times \N^k\subset \mathds{H}^{N}\times \R^k$ and  the infinitesimal generator for $u^N$ has the form
\[
\begin{array}{l}
\displaystyle \mathcal{A}^N\varphi(u_C,u_D)\hspace{1cm}\vspace{0.2cm}\\  
\displaystyle \hspace{0.2cm}=\sum_{j=1}^N\left\lbrace\sum_{r\in\mathfrak{R}_C}\left[\varphi\left(u_C+
\frac{\gamma_r^C}{\mu}\mathds{1}_j,u_D\right)-\varphi(u_C,u_D)\right]\mu\lambda_r\left(u_j^{C}\right)\right.\\ 

\displaystyle \hspace{1.7cm}\left.+\sum_{r\in S_1} \left[\varphi\left(u_C+\frac{\gamma_r^C}{\mu}\mathds{1}_j,u_D\right)-\varphi(u_C,u_D)\right]\mu\lambda_r\left(u_j^{C},u_{\ell_j}^{D}\right)\right\}\\ 

\displaystyle \hspace{0.3cm}+\sum_{j=1}^N\left\lbrace\left[\varphi\left(u_C+\frac{\mathds{1}_{j-1}-\mathds{1}_j}{\mu},u_D\right)+\varphi\left(u_C+\frac{\mathds{1}_{j+1}-\mathds{1}_j}{\mu},u_D\right)-2\varphi(u_C,u_D)\right]\mu N^2u_j^C\right\}
\end{array}
\]
\begin{equation}\label{generator_of_the_JMP}
\begin{array}{l}
\displaystyle \hspace{0.3cm}+\sum_{\ell=1}^k\left\{\sum_{r\in\mathfrak{R}_{DC}\backslash S_1}\left[\varphi\left(u_C+\frac{\gamma_{r}^{C,N}}{\mu},u_D+\gamma_{r}^{D}\1^\ell\right)-\varphi(u_C,u_D)\right]\lambda_r\left(\sum_{j=N_{\ell-1}+1}^{N_\ell}a_{j,N}^ru_j^{C},u^{D}_\ell\right)\right. \\

\displaystyle \hspace{2cm}\left.+\sum_{r\in\mathfrak{R}_D}\left[\varphi\left(u_C,u_D+\gamma_{r}^{D}\1^\ell\right)-\varphi(u_C,u_D)\right]\lambda_r\left(u_\ell^{D}\right)\right\rbrace
\end{array}
\end{equation}

\noindent on the domain $\displaystyle C_b\left(\mathds{H}^{N}\times\R^k\right)$. It can be extended to an $\tilde{\cA}^N$ on $\disp C_b\big(L^2(I)\times \R^k\big)$ by 
\begin{equation}\label{extended_generator_of_the_JMP}
\disp \tilde{\cA}^N\varphi(u_C,u_D)=\cA^N\varphi(P_Nu_C,u_D).
\end{equation}
\noindent Again, we omit the tildes: $\tilde{\cA}^N=\cA^N$. Also we write $\tilde{P}_Nu:=\big(P_Nu_C,\sum_{\ell=1}^ku_{\ell}^D\1^\ell)$ for $u=(u_C,u_D)\in L^2(I)\times \R^k$.\\


\subsection{Convergence tools}

\subsubsection{Formal limit of the generator}\label{formal_limit_of_the_generator2}
\label{s2.4.1}

We seek for the asymptotic behavior of the process $u^N$ presented above, as $N,\mu\rightarrow\infty$. We first introduce some mathematical tools.

For $f,g\in\mathds{H}^N$, the $L^2$ inner product reads $\left\langle f,g\right\rangle_2=N^{-1}\sum_{j=1}^Nf_jg_j$, and the supremum norm is given by $\Vert f\Vert_\infty=\max_{1\leq j\leq N}|f_j|$. 

\begin{proposition}\label{Piecewise_Approximation_via_Projection}\hspace{1cm}\vspace{0.1cm}\par 

\textbf{(i)} $\displaystyle \left(\mathds{H}^N,\langle\cdot,\cdot\rangle_2\right)$ is a finite dimensional Hilbert space with $\displaystyle \left\{\sqrt{N}\mathds{1}_j, 1\leq j\leq N\right\}$ as an orthonormal basis.\par 
\textbf{(ii)} $\displaystyle \lim_{N\rightarrow\infty}\Vert P_Nf-f\Vert_2\longrightarrow 0$ for $f\in L^2$, 
\par 
\textbf{(iii)} $\displaystyle \lim_{N\rightarrow\infty}\Vert P_Nf-f\Vert_\infty\longrightarrow 0$, for $f\in C(I)$ .

\end{proposition}

Henceforth, $\H^N\times\R^k$ is endowed with the norm \[ \Vert(f,g)\Vert_{\infty,\infty}:= \Vert f\Vert_\infty+\Vert g\Vert_\infty,\hspace{0.5cm}\text{for}\hspace{0.2cm}f\in \H^N,\; g\in \R^k. \]

\noindent \textbf{\textit{The discrete Laplace.}} For $f\in D(I)$ and for $x\in[0,1]$, we set \[ \nabla_N^+ f(x):=N\left[f\left(x+\frac{1}{N}\right)-f(x)\right]\hspace{0.5cm}\text{and}\hspace{0.5cm}\nabla_N^- f(x):=N\left[f(x)-f\left(x-\frac{1}{N}\right)\right]. \] Then, we define the discrete Laplace on $D(I)$ by\\

$\displaystyle \hspace{1cm} \Delta_Nf(x):=\nabla_N^+\nabla_N^-f(x)$\par 
$\displaystyle \hspace{2.7cm}=\nabla_N^-\nabla_N^+f(x)=N^2\left[f\left(x-\frac{1}{N}\right)-2f(x)+f\left(x+\frac{1}{N}\right)\right]$.\\

\noindent If $f\in\H^N$ in particular, then \[ \Delta_Nf(x)=\sum_{j=1}^N\big[N^2(f_{j-1}-2f_j+f_{j+1})\big]\1_j(x). \] 

From the spectral analysis of $\Delta_N$, it is well known that, if $N$ is an odd integer, letting $0\leq m\leq N-1$ with $m$ even, letting $\varphi_{0,N}\equiv 1$, 
$\displaystyle \varphi_{m,N}(x)=\sqrt{2}cos\big(\pi mjN^{-1}\big)$ and $\displaystyle \psi_{m,N}(x)=\sqrt{2}sin\big(\pi mjN^{-1}\big)$ for $\displaystyle x\in I_j$, then, $\displaystyle \big\{\varphi_{m,N},\psi_{m,N}\big\}$ are eigenfunctions of $\Delta_N$ with eigenvalues given by $\displaystyle -\beta_{m,N}=-2N^2\left(1-cos\big(\pi mN^{-1}\big)\right)\leq 0$. If $N$ is even, we need the additional eigenfunction $\displaystyle \varphi_{N,N}=cos(\pi j)$ for $x\in I_j$. The following (classical) properties are derived from \cite{Blount1987}, Lemma 2.12 p.12, \cite{blount1992}, Lemma 4.2 for the parts (i)-(v), from \cite{Kato1966}, Chapter 9, Section 3 for the part (vi), and from \cite{Henry1981}, Chapter 1, Sections 1.3 and 1.4 for the parts (viii)-(ix). Consider $f,g\in\mathds{H}^N$ and let $\displaystyle T_N(t)=e^{\Delta_Nt}$ denote the semigroup on $\mathds{H}^N$ generated by $\Delta_N$.\vspace{0.1cm}\par 

\begin{proposition}\label{Discrete_Laplacian_properties}

\noindent\hspace{0.2cm} \textbf{(i)} The family $\{\varphi_{m,N},\psi_{m,N}\}$ forms an orthonormal basis of $\displaystyle \left(\mathds{H}^N,\langle\cdot,\cdot\rangle_2\right)$.\vspace{0.1cm}\par 

\noindent\hspace{0.2cm} \textbf{(ii)} $\disp T_N(t)f=\sum_{m}e^{-\beta_{m,N}t}\left(\langle f,\varphi_{m,N}\rangle_2\varphi_{m,N}+\langle f,\psi_{m,N}\rangle_2\psi_{m,N}\right)$.\vspace{0cm}\par 

\noindent\hspace{0.2cm} \textbf{(iii)} $\langle \nabla_N^+f,g\rangle_2=\langle f,\nabla_N^-g\rangle_2$\hspace{0.1cm} and \hspace{0.1cm}$T_N(t)\Delta_Nf=\Delta_NT_N(t)f$. \vspace{0.1cm}\par 

\noindent\hspace{0.2cm} \textbf{(iv)} $\Delta_N$ and $T_N(t)$ are self-adjoint on $\displaystyle \left(\mathds{H}^N,\langle\cdot,\cdot\rangle_2\right)$.\vspace{0.1cm}\par 

\noindent\hspace{0.2cm} \textbf{(v)} $T_N(t)$ is a positive contraction semigroup on both $\displaystyle \left(\mathds{H}^N,\langle\cdot,\cdot\rangle_2\right)$ and $\displaystyle \left(\mathds{H}^N,\Vert\cdot\Vert_\infty\right)$.\vspace{0.1cm}\par 

\noindent\hspace{0.2cm} \textbf{(vi)} The projection $P_N$ commutes with $\Delta_N$, and for all $f\in C^2(I)$, $t\geq0$,
\[ \left\{
\begin{array}{l}
\disp \Vert \Delta_Nf-\Delta f\Vert_\infty\longrightarrow0\hspace{0.3cm}\text{as}\hspace{0.2cm}N\rightarrow\infty,\vspace{0.1cm}\\
\disp \Vert \Delta_NP_Nf-\Delta f\Vert_\infty\longrightarrow0\hspace{0.3cm}\text{as}\hspace{0.2cm}N\rightarrow\infty,\vspace{0.1cm}\\
\disp \Vert T_N(t)P_Nf-T(t)f\Vert_\infty\longrightarrow0\hspace{0.3cm}\text{as}\hspace{0.2cm}N\rightarrow\infty.
\end{array} 
\right. \]

\noindent\hspace{0.2cm} \textbf{(vii)} Let $f_N,f\in C^2(I)$.  \[ \text{If }\hspace{0.2cm} \lim_{N\rightarrow\infty}\Vert f_N-f\Vert_{C^2(I)}\rightarrow 0,\hspace{0.4cm} \text{then }\hspace{0.5cm}\lim_{N\rightarrow\infty} \Vert\Delta_Nf_N-\Delta f\Vert_\infty\rightarrow0. \]

\noindent\hspace{0.2cm} \textbf{(viii)} For all $\eta\geq 0$, there exist a constant $c_1=c_1(\eta)>0$ independent of $N$, such that \[ \left\Vert T_N(t)\big(-\Delta_N\big)^{\eta}\right \Vert_{L^\infty(I)\rightarrow L^\infty(I)}\leq c_1t^{-\eta}\hspace{0.4cm}\text{for}\hspace{0.5cm} t>0. \]

\noindent\hspace{0.2cm} \textbf{(ix)} For all $0< \eta\leq 1$, there exists $c_2=c_2(\eta)>0$ independent of $N$, such that \[ \left\Vert \left(T_N(t)-I_d\right)\big(-\Delta_N\big)^{-\eta}\right\Vert_{L^\infty(I)\rightarrow L^\infty(I)}\leq c_2t^\eta \hspace{0.4cm}\text{for}\hspace{0.5cm}t\geq0. \]

\end{proposition}

We wish to identify the limit as $N,\mu\rightarrow\infty$, of our sequence $u^N$ of Markov processes with generators $\cA$ given by $(\ref{generator_of_the_JMP})$, we compute a formal limit of the sequence of the extended generators defined by $(\ref{generator_of_the_JMP})$, for test functions $\disp \varphi\in C_b^{2,0}(L^2(I)\times \R^k)$. As already mentioned the rigorous proof below is done for a stronger topology than $L^2(I)$, but this latter space being Hilbert the formal argument is clearer. 

We first proceed to a Taylor expansion at order $2$ of $\varphi$ for the part of the generator corresponding to the fast dynamic, that is the first three terms in the right hand side of (\ref{generator_of_the_JMP}). We identify the differential with the corresponding gradient, the debit functions appear naturally. With the notation $u=(u_C,u_D)$, we have: \\

\noindent \textbf{\textit{First order terms.}}\vspace{0.1cm}\par 

$\disp\hspace{0cm} T_{fast}^{(1)}(N)= \left\langle D^{1,0}\varphi\big(\tilde{P}_Nu\big),\sum_{j=1}^N\left(\sum_{r\in\fR_C}\gamma_r^C\lambda_r\big(u_j^C\big)+\sum_{r\in S_1}\gamma_r^C\lambda_r\big(u_j^C,u_{\ell_j}^D\big)\right)\1_j\right\rangle_2$,\par  
$\disp \hspace{2.5cm}+\left\langle D^{1,0}\varphi\big(\tilde{P}_Nu\big),\sum_{j=1}^NN^2\left(\mathds{1}_{j-1}-2\mathds{1}_{j}+\mathds{1}_{j+1}\right)u_j^C\right\rangle_2$\par 

$\disp \hspace{1.5cm}=\left\langle D^{1,0}\varphi\big(\tilde{P}_Nu\big),F\big(\tilde{P}_Nu\big)+\Delta_NP_Nu_C\right\rangle_2$.\\

\noindent Here $\Delta_N$ is the discrete Laplace introduced above and corresponds to the debit function corresponding to the diffusions. Its expression is obtained using a change of index and periodicity. The function $F$ is the debit related to fast onsite reactions. It also maps $D(I)\times \R^k$ on $D(I)$, and is given by 
\begin{equation}\label{debit_function_of_fast_onsite_reactions}
\disp F\left(u_C,u_D\right)=\sum_{r\in\mathfrak{R}_C}\gamma_r^C\lambda_r\left(u_C\right)
+\sum_{r\in S_1}\gamma_r^C\lambda_r\left(u_C,u_D\right).
\end{equation}
These allow to define the debit of the whole fast dynamic $\disp (u_C,u_D)\mapsto \Delta_NP_Nu_C+F\big(\tilde{P}_Nu\big)$ on $D(I)\times \R^k$. If $F$ is continuous, then formally for $u_C\in C^2(I)$,
\begin{equation}\label{formal_limit_of_the_fast_part_of_the_generator}
\disp T_{fast}^{(1)}(N)\longrightarrow \big\langle D^{1,0}\varphi(u_C,u_D), \Delta u_C+F(u_C,u_D)\big\rangle_2=:\cA_{fast}^\infty
\end{equation}
as $N\rightarrow\infty$, thanks to $\cref{Discrete_Laplacian_properties}$ (vi) and $\cref{Piecewise_Approximation_via_Projection}$ (iii). \\

\noindent \textbf{\textit{Second order terms.}} We give heuristics indicating that they converge to zero in the $L^2(I)$ topology. For these heuristics, we assume that the concentration of the species $C$ is bounded (uniformly in $N$) on every site. Say $0\leq u_C(x)\leq\rho$ for all $x\in I$ for some $\rho>0$. Then in particular, the rates of fast onsite reactions are bounded too, if these latter are continuous. Let us introduce $e_j:=\sqrt{N}\1_j$, for $1\leq j\leq N$. Since the fonctions $\disp (\1_j)_{1\leq j\leq N}$ are pairwise orthogonal in $L^2(I)$ and of norm $\disp \Vert \1_j\Vert_2=N^{-1/2}$, the family $\big\{e_j,1\leq j\leq N\big\}$ forms an orthonormal basis of $\disp \big(\H^N,\Vert\cdot\Vert_2\big)$. Also, $c$ denotes a generic  constant and $\Vert D^{2,0}\varphi\Vert_\infty$ denotes a uniform bound on the norm of  $D^{2,0} \varphi(u)$ as 
a bilinear operator on $L^2(I)$ for $u\in L^2(I)\times \R^k$.

For fast onsite reactions, we have the following term for the second order term in the Taylor expansion\\

$\disp\hspace{0.5cm} T_{fastOn}^{(2)}(N)=\frac{1}{2}\sum_{j=1}^N\left(\sum_{r\in \fR_C}D^{2,0}\varphi\big(\tilde{P}_Nu\big)\cdot\left\langle\frac{\gamma_r^C}{\mu}\mathds{1}_j,\frac{\gamma_r^C}{\mu}\mathds{1}_j\right\rangle_2\mu\lambda_r\left(u_j^C\right)\right.$\par

$\disp\hspace{4.3cm}\left.+\sum_{r\in S_1}D^{2,0}\varphi\big(\tilde{P}_Nu\big)\cdot\left\langle\frac{\gamma_r^C}{\mu}\mathds{1}_j,\frac{\gamma_r^C}{\mu}\mathds{1}_j\right\rangle_2\mu\lambda_r\left(u_j^C,u_{\ell_j}^D\right)\right)$

$\disp \hspace{2.5cm}\leq c\bar{\gamma}_C^2\bar{\lambda}_C\left\Vert D^{2,0}\varphi\right\Vert_\infty\frac{1}{N\mu}\sum_{j=1}^N\big\Vert e_j\big\Vert_2^2\hspace{0.2cm}\leq \hspace{0.2cm} \frac{c}{\mu}\longrightarrow0$.\\

\noindent Here, $\bar{\gamma}_C$ is an upper bound for the amplitudes of fast onsite reactions and $\bar{\lambda}_C$ is an upper bound for the rates of fast onsite reactions. 

We now treat the diffusion part. From $\disp \Vert\1_{j+1}-\1_j\Vert_2^2=2N^{-1}=\Vert\1_{j-1}-\1_j\Vert_2^2$,\\

$\disp \hspace{0.5cm} T_{fastDiff}^{(2)}(N)=\frac{1}{2}\sum_{j=1}^N\left[D^{2,0}\varphi\big(\tilde{P}_Nu\big)\cdot\left\langle\frac{\mathds{1}_{j-1}-\mathds{1}_{j}}{\mu},\frac{\mathds{1}_{j-1}-\mathds{1}_{j}}{\mu}\right\rangle_2\right]\mu N^2u_j^C$\par 

$\disp \hspace{3cm}\leq c\rho\left\Vert D^{2,0}\varphi\right\Vert_\infty\frac{N^2}{\mu}\sum_{j=1}^N\Vert\1_{j-1}-\1_j\Vert_2^2$\par 

$\disp \hspace{3cm}\leq c\mu^{-1}N^2\longrightarrow0$,\\

\noindent if $\disp\mu^{-1}N^2\rightarrow0$ as $N,\mu\rightarrow\infty$. 
The other term related to diffusions 
converges to zero, using the same argument.  Below, we will see that the condition $\mu^{-1}N^2\rightarrow\infty$ can be weakened.\\

Let us then treat the slow dynamic part. With our choice of $b_{j,N}^r$, $\gamma_{r,\ell}^{C,N}$ are clearly bounded and $\frac{\gamma_{r,\ell}^{C,N}}\mu\to 0$. Also, since $a^r$ is continuous, it is not difficult to
see that if $u_C^{N}$ and $u_\ell^{N,D}$ converge to $u_C$ and  $u_{D,\ell}$ respectively, then $\lambda_r\left(\sum_{j=N_{\ell-1}+1}^{N_\ell}a_{j,N}^r u_j^{N,C},u^{N,D}_\ell\right)$ converges to $\lambda_r\left(\int_{J_\ell}a^r(x)u_C(x)dx,u_{D,\ell}\right)$. Since here, we work at a formal level, we do not need to precise in which sense the convergences hold.

We deduce the limit generator:
\begin{equation}\label{limiting_generator}
\begin{array}{ll}
\displaystyle \mathcal{A}^\infty\varphi(u_C,u_D) &=\left\langle D_{u_C}\varphi(u_C,u_D),\Delta u_C+F(u_C,u_D)\right\rangle_2\\
&\displaystyle +\sum_{\ell=1}^k\left\{\sum_{r\in\mathfrak{R}_{DC}\backslash S_1}\left[\varphi\left(u_C,u_D+
\gamma_r^D \1^\ell\right)-\varphi(u_C,u_D)\right]\lambda_r\left(\int_{J\ell}a^r(x)u_C(x)dx,u_{D,\ell}\right)\right.\vspace{0.1cm}\\
&\displaystyle \hspace{1.4cm}\left.+\sum_{r\in\mathfrak{R}_D}\left[\varphi\left(u_C,u_D+\gamma_r^D\1_\ell\right)-\varphi(u_C,u_D)\right]\lambda_r\left(u_{D,\ell}\right)\right\rbrace
\end{array}\\
\end{equation}

\noindent As seen above, this is the generator of a continuous PDMP.

\subsubsection{The problem at the limit}

From $\cref{infinite_dimensional_PDMPs}$, we know that $(\ref{limiting_generator})$ is the (extended) generator of a continuous PDMP $v:=\{v(t)=(v_C(t),v_D(t)),t\geq0\}$ as given by $(\ref{extended_generator_of_a_PDMP})$. The evolution of its discrete component $v_D(t)$ is governed by the slow dynamic part in the limiting generator. It is a jump process, and some of its jump rates depend on both components. The possible values of $v_D(t)$ are of the form 
$$
v_D=v_D(0)+\sum_{m=1}^M\sum_{r\in (\mathfrak{R}_{DC}\backslash S_1)\cup \mathfrak{R}_{D}}\sum_{\ell=1}^k \alpha_{\ell,r}^m\gamma_r^D\1^\ell,
$$
where the coefficients $\alpha_{\ell,r}^m$ are integers. This describes a countable set denoted by $K$.

The continuous component $v_C(t)$ has values in $C(I)$, and has continuous trajectories. More precisely, for any value $\nu\in K$ of $v_D(t)$ ---fixed between two consecutive jumps of this latter---, $v_C(t)$ evolves according to
\begin{equation}\label{abstract_evolution_equation_of_the_limiting_PDMP}
\displaystyle \frac{d}{dt}v_C(t)=\Delta v_C(t)+F\left(v_C(t),\nu\right), 
\end{equation}
the reaction-diffusion equation with a parameter $\nu$. It is supplemented with periodic boundary conditions. We consider the Laplace operator $\Delta$ as an operator in the $C$-theory framework, w.r.t. the terminology used in \cite{Cazenave1998} (Chapter 2, Section 6). Dirichlet boundary conditions are considered
there, the case of periodic boundary conditions considered here is similar. 
We use the same symbol $\Delta$ for this operator on the domain:
 \[ \big\{u\in C(I):\Delta u\in C(I)\big\}. \] This latter domain is dense, and the Laplace is also $m$-dissipative. 
From Hille-Yosida-Phillips theorem, it is well known that this operator generates  a strongly continuous semigroup of contraction on $C(I)$, denoted by $\{T(t)\}=\left\{T(t)=e^{\Delta t}, t\geq0\right\}$.

\begin{assumption}\label{well-posedness_assumpt_on_the_deb_funct_of_fast_onsite_react}
\hspace{1cm}\par  
(i) For all $y=(y_1,y_2)\in\R^2$, $F(y)\geq 0$ if $y_1=0$.\par 
(ii) There exists $\rho_1>0$ such that $F(y)<0$ for all $y\in\R^2$ satisfying $|y_1|>\rho_1$. 
\end{assumption}
\noindent Under $\cref{well-posedness_assumpt_on_the_deb_funct_of_fast_onsite_react}$, and given any initial condition $v_C(0)=v_{C,0}\in C(I)$ such that $0\leq v_{C,0}(x)\leq \rho_2$, for all $x\in I$, for some $\rho_2>0$, there exists a unique global mild 
solution $v_C$ to $(\ref{abstract_evolution_equation_of_the_limiting_PDMP})$, satisfying $\disp v_C\in C\big(\R_+;C(I)\big)$, and $0\leq v_C(t,x)<\rho=(\rho_1\vee\rho_2)+1$: for all $t\geq 0$, $x\in I$. We use the notation $v_C(t)=\phi_\nu(t)v_{C,0}$ for the associated flow.\par 
By mild solution, we mean for all $t\ge 0$:
$$
v_C(t)=T(t)v_{C,0}+\int_0^t T(t-s)F(v_C(s),\nu)ds.
$$
Recall that $F$ is polynomial. 
Note that if $v_{C,0}\in W^{2,\infty}(I)$ then $v_C$ belongs to $C([0,T];W^{2,\infty}(I))$. This can be seen from 
standard properties of the heat kernel. Moreover, there exists a constant depending on $\rho$ and $F$ such that:
\begin{equation}\label{e2.19}
\|v_C(t)\|_{W^{2,\infty}}\le c(\rho,F,\|v_{C,0}\|_{W^{2,\infty}}),\; t\ge 0.
\end{equation}

From the above discussion, we consider that the "limiting" PDMP $v=(v_C,v_D)$ starts at $v_0=(v_{C,0},v_{D,0})\in C(I)\times \R^k$, has values in  
$ E:=C(I)\times K$. Below, we denote by $\|\cdot\|_\infty$ the norm of  $C(I)\times \R^k$. It has the following characteristics: the \textit{linear and nonlinear operators} are respectively the Laplace operator $\Delta$ and the debit of fast onsite reactions $F$ given by $(\ref{debit_function_of_fast_onsite_reactions})$, the \textit{transition rate} $\Lambda$, defined by 
\begin{equation}\label{global_jump_rate_slow_reaction}
\displaystyle\hspace{0cm} \Lambda(u_C,u_D):=\sum_{\ell=1}^k\left\{\sum_{r\in\mathfrak{R}_{DC}\backslash S_1}\lambda_r\left(\int_{J_\ell}a^r(y)u_C(y)dy,u_{D,\ell}\right)+\sum_{r\in\mathfrak{R}_{D}}\lambda_r\left(u_{D,\ell}\right)\right\},
\end{equation} 
and, the \textit{transition measure} $Q$, defined by 
\begin{equation}\label{transition_measure}
\begin{array}{l}
\displaystyle \int_{K}\varphi(u_C,\nu)Q\big(d\nu;(u_C,u_D)\big)\vspace{0.1cm}\\

\displaystyle \hspace{0.5cm}=\frac{1}{\Lambda(u_C,u_D)}\hspace{0.1cm}\sum_{\ell=1}^k\left\{\sum_{r\in\mathfrak{R}_{DC}\backslash S_1}\varphi\left(u_C,u_D+\gamma_r^D\1^\ell\right)\lambda_r\left(\int_{J_\ell}a^r(y)u_C(y)dy,u_{D,\ell}\right)\right.\vspace{0.1cm}\\

\displaystyle \hspace{4cm}\left.+\sum_{r\in\mathfrak{R}_D}\varphi\left(u_C,u_D+\gamma_r^D\1^\ell\right)\lambda_r\left(u_{D,\ell}\right)\right\rbrace.
\end{array}
\end{equation}
We assume that \cref{non_explosion_of_the_PDMP} holds. Then, $v_D$ is a c\`adl\`ag process which is almost surely in $D([0,T],\R^k)$. 

Below, during the proof, we see that we may assume that all reaction rates are bounded. Then, arguing as in \cite{Davis1993} or \cite{Riedler2011bis}, we prove that the law of $(v_C,v_D)$ solves the martingale problem in  the sense that for all $\varphi\in \mathcal E$:
$$
\varphi(v_C(t),v_D(t))-\varphi(v_{C,0},v_{D,0})-\int_0^t \mathcal{A}^\infty\varphi(v_C(s),v_D(s))ds
$$
is a martingale.
The set of test functions $\mathcal E$ is the set of functions $\varphi$ on $E=C(I)\times \R^k$  which are bounded, 
continuous with respect to the second variable, 
differentiable with bounded differential with respect to the first variable and such that
the differential can be extended continuously to $H^{-1}(I)$. Recall that this is the dual space of $H^1_0(I)$, and there exists $c_\varphi$ such that 
$$
|D_{u_C}\varphi(u_C,u_D)\cdot h| = |\langle D_{u_C}\varphi(u_C,u_D), h\rangle_2| \le c_\varphi \|h\|_{H^{-1}(I)}
$$
for all $(u_C,u_D)\in C(I)\times \R^k$, $h\in L^\infty(I)$. We denote by $|\varphi|_{\cE}$ the smallest constant $c_\varphi$ such that the above inequality holds and
$$
\|\varphi\|_\cE= \|\varphi\|_{C^{1,0}_b(C(I)\times \R^k)}+|\varphi|_\cE.
$$

The forthcoming result states that, under the boundedness assumption, the martingale problem for $\cA^\infty$ is well-posed.

\begin{theorem}\label{well_posedness_of_the_martingale_problem_for_inf_dim_PDMP} Assume that the reaction rates are bounded, the (law of the) PDMP $v$  is the unique solution of the martingale problem associated with   $\cA^\infty$ on  $\mathcal E$.
\end{theorem}

The proof is posponed to section $\ref{proof_of_well_posedness_of_the_martingale_problem_for_inf_dim_PDMP}$. It generalizes the proof of Theorem 2.5 in \cite{Arnaud2012}.


\section{A Piecewise deterministic limit}\label{s3}

The following is our main result.

\begin{theorem}\label{hybrid_simplification}\textbf{\texttt{(Hybrid simplification)}}\vspace{0.1cm}\par

\noindent Let $\displaystyle u^N=\left(u_C^N,u_D^N\right)$ define a sequence of Markov processes as above, starting at $\displaystyle u^N(0)=\left(u_C^N(0),u_D^N(0)\right)\in \H^N\times K$,
with the corresponding sequence of infinitesimal generators $\mathcal{A}^N$, defined by $(\ref{generator_of_the_JMP})$. Assume the jump rates $\lambda_r$, $r\in\fR$ are polynomial, and the rates $\lambda_r$, $r\in\fR_C\cup S_1$ are such that $F$ satisfies $\cref{well-posedness_assumpt_on_the_deb_funct_of_fast_onsite_react}$. \par 

Now, consider the $E-$valued PDMP $v=(v_C,v_D)$  started from $v_0=(v_0^C,v_0^D)\in C^2(I)\times K$, whose characteristics are $\left(\Delta,F,\lambda^D,Q\right)$. Assume \cref{non_explosion_of_the_PDMP} holds so that the PDMP is standard.\par 

In addition, assume that: (i) $N,\mu\rightarrow\infty$ such that $\displaystyle \mu^{-1}\log N\longrightarrow0$.\par 
\hspace{4.5cm} (ii) $\displaystyle \left\Vert u^N(0)-v_0\right\Vert_{\infty}\longrightarrow 0$ in probability. \\
Then for all $T>0$, \[ u^N\Longrightarrow v \hspace{0.5cm}\text{ in }\hspace{0.3cm}D\left([0,T];D(I)\times K\right), \]
where $\Longrightarrow$ means convergence in law.
\end{theorem}\vspace{0.4cm}\par

\noindent \textbf{\textit{Proof}.} 
Relying on the proof of Theorem 3.1 in \cite{Arnaud2012}, we proceed as follows. {We first assume that reaction rates} corresponding to onsite reactions are bounded, as well as their derivatives with respect to their first variable. In this case, 
we prove tightness for the sequence of Markov processes and identify its limit through the associated sequence of martingale problems. 

At the end of the proof, we relax the additional hypothesis of boundedness and conclude by a truncation argument.\par 

In the whole proof, $T>0$ is fixed and {every limit is taken} as $N\rightarrow\infty$ if there is no further precision. \\


\subsection{Tightness for bounded reaction rates}

We want to show that the family of processes $\displaystyle \left\{u^N=\left(u_C^N,u_D^N\right)\right\}_N$ is tight in the Skorohod's space $\displaystyle D([0,T],D(I)\times
\R^k)$. This has to be understood in the sense that the family $\displaystyle \left\{P_N\right\}_N$ of their laws ---where $P_N$ is the law of $u^N$ for all $N$--- is tight in 
the space $\displaystyle \mathcal{P}\big(D([0,T],D(I)\times \R^k)\big)$ of probability measures (see \cite{Billingsley1999}, page 8). We treat both components one after another, starting with the continuous one.

We assume until section \ref{s3.3} that all reaction rates are bounded.

\subsubsection{Tightness for the continuous component} \label{s3.1.1}

We first prove that $\displaystyle \left\{u_C^N\right\}_N$ is $C$-tight. That is (see \cite{Jacod1987}, Chapter 6, Section 3, Definition 3.25, page 351), 
the family is tight in $D([0,T],D(I))$, and each converging subsequence converges to a limit whose law has its support in $C([0,T],D(I))$, {\it i.e.} the limiting law charges only continuous trajectories. 
In fact we even prove that each limit point has support in $C([0,T],C(I))$. 

From Proposition 2.1 of \cite{Kurtz1971}, we know that $u_C^N$ is solution to the SDE
\[ u_C^N(t):=u_C^N(0)+\int_0^t\Delta_N u_C^N(s)+F\big(u^N(s)\big) ds+Z_C^N(t), \] where $Z_C^N$ is a $\displaystyle \mathds{H}^N$-valued martingale on $[0,T]$. 

Duhamel's formula yields 
$$
\displaystyle u_C^N(t)=T_N(t)u_C^N(0)+\int_0^tT_N(t-s)F\big(u^N(s)\big)ds+Y_C^N(t),
$$
where $\displaystyle Y_C^N(t)=\int_0^tT_N(t-s)dZ_C^N(s)$ is a stochastic convolution. 

As in \cite{blount1992} and \cite{DebusscheNankep2017}, we need to introduce a stopping time to control $Y_C^N$. Define $v^N_C$ by
\begin{equation}\label{e3.2bis}
v^N_C(t)=T_N(t)u_C^N(0)+\int_0^tT_N(t-s)F\big(v^N_C(s),u_D^N(s)\big)ds.
\end{equation}
It is classical to prove that $v^N_C$ exists on $\R^+$ and since $F$ is bounded and $(T_N(t))_t$ is a contraction semigroup, we have:
\begin{equation}\label{e3.2}
\|v^N_C(t)\|_\infty\le \|u^N_C(0)\|_\infty +T M_F , \quad t\in [0,T],
\end{equation}
where $M_F$ is a bound on $F$. Recall that $F$ is defined in \eqref{debit_function_of_fast_onsite_reactions} and all reaction rates are assumed to be bounded 
in this part of the proof.

We also introduce 
\begin{equation}\label{SDE_for_the_truncated_continous_Markov_process}
\bar u_C^N(t)=T_N(t)u_C^N(0)+\int_0^tT_N(t-s)F\big(\bar u^N_C(s),u^N_D(s)\big)ds+Y_C^N(t)\1_{t\le \tau_N}
\end{equation}
where the stopping time $\tau_N$ is defined by 
$$
\tau_N=\inf\{t\ge 0, \; \|u^N_C(t)-v^N_C(t)\|_\infty\ge 1\}.
$$
Using the same arguments as in  \cite{blount1992} and \cite{DebusscheNankep2017}, we prove
 \[ \sup_{t\in[0,T\wedge \tau_N ]}\left\Vert Y_C^N(t)\right\Vert_\infty\longrightarrow0\hspace{0.5cm}\text{in probability}, \] 
 under the assumptions of \cref{hybrid_simplification}. Using a Gronwall type argument, we deduce:
 \begin{equation}\label{e3.4}
 \sup_{t\in [0,T]} \|\bar u^N_C(t)-v^N_C(t)\|_\infty \longrightarrow0\hspace{0.5cm}\text{in probability}.
 \end{equation}
Then, we write 
$$
\P\left( \sup_{t\in [0,T]} \|\bar u^N_C(t)-u^N_C(t)\|_\infty >\epsilon \right) \le \P\left( \sup_{t\in [0,T]} \|\bar u^N_C(t)-v^N_C(t)\|_\infty \ge 1 \right)\to 0, 
$$
when $N\to \infty$ and, by Lemma 3.31, Section 3, Chapter 6 of \cite{Jacod1987},  tightness of $\{ u_C^N\}_N$ follows form tightness of $\{{\bar u}_C^N\}_N$

Under our assumptions:
$$
T_N(t)u_C^N(0) \to T(t) v_C(0), \mbox{ in } C(I)
$$
so that by \eqref{e3.4}, it it suffices to prove tightness of 
%
%
%
%
%
$$\disp b^N(t):=\int_0^tT_N(t-s)F\big(\bar u^N(s),u_D^N(s)\big)ds.$$

%
We use Arzel\`a-Ascoli theorem to get relative compactness and derive tightness from Prohorov theorem, since $\displaystyle C\big([0,T],D(I)\big)$ is separable and complete.\\

\noindent \textbf{\textit{Equicontinuity (in time).}} Fix $t_1,t_2\in[0,T]$. 
%
Let $\displaystyle I_d$ be the identity operator on $\displaystyle D(I)$, and let $0<\eta< 1$ be a real number. We have\\

$\disp \hspace{1cm} b^N(t_2)-b^N(t_1) = \int_0^{t_1}\left(T_N(t_2-t_1)-I_d\right)T_N(t_1-s)F\big(\bar u^N_C(s),u^N_D(s)\big)ds$\par 
$\disp \hspace{5cm} +\int_{t_1}^{t_2}T_N(t_2)F\big(\bar u^N_C(s),u^N_D(s)\big)ds$\\

From \cref{Discrete_Laplacian_properties} (iii), (viii) and (ix), and the boundedness of $F$:\\

$\displaystyle \left\Vert b^{N}(t_2)-b^{N}(t_1)\right\Vert_{\infty}$\par 

$\displaystyle \hspace{1cm} \leq\int_0^{t_1}\left\Vert\big(T_N(t_2-t_1)-I_d\big)\big[\big(-\Delta_N\big)^{-\eta}\big]\big[\big(-\Delta_N\big)^{\eta}\big]T_N(t_1-s)F
\big(\bar u^N_C(s),u^N_D(s)\big)\right\Vert_\infty ds$\vspace{0.1cm}\par 
$\displaystyle \hspace{3cm}+\int_{t_1}^{t_2}\left\Vert T_N(t_2)F\big(\bar u^N_C(s),u^N_D(s)\big)\right\Vert_\infty ds$\par 

$\displaystyle \hspace{1cm} \leq c_2|t_2-t_1|^\eta M_F\int_0^{t_1}c_1(t_1-s)^{-\eta}ds+|t_2-t_1|M_F$,\hspace{1cm} $\forall 0<\eta< 1$,\par 

$\displaystyle \hspace{1cm}\leq c|t_2-t_1|^\eta M_F\left(\frac{t_1^{1-\eta}}{1-\eta}+|t_2-t_1|^{1-\eta}\right)$, \hspace{1cm}$c=c_1\vee c_2=\max(c_1,c_2)$\\

$\displaystyle \hspace{1cm}\leq c|t_2-t_1|^\eta$, \hspace{1cm} $c=(M_F,c_1,c_2,T)$,\\

\noindent and the family $\displaystyle \left\{b^N\right\}_N$ is uniformly equicontinuous. \\


\noindent \textbf{\textit{Compactness (in space).}} Fix $t\in[0,T]$. 

Using the contraction $T_N(t)$ property and boundedness of $F$, we get 
\[ \sup_{t\in[0,T]}\Vert b^N(t)\Vert_\infty\leq TM_F=:c_5. \]

Let $i_0$ be such that  $(\nabla_N^+v)_{i_0}= \Vert\nabla_N^+v\Vert_\infty$ and write:
$$
(\nabla_N^+v)_{i_0}= (\nabla_N^+v)_{i_0-1}+N (\Delta_N v)_{i_0}
$$

By \cite{beale2009}, Theorem 3.1, we know that for all $t>0$, $\Vert\nabla_N^+T_N(t)\Vert_{L^\infty(I)\rightarrow L^\infty(I)}\leq ct^{-1/2}$ for some constant $c$ independent of $N$. Since $F$ is bounded, we have
\begin{align*}
\Vert\nabla_N^+b^N(t)\Vert_\infty & \leq \int_0^t\Vert \nabla_N^+T_N(t-s)F(\bar u^N_C(s),u^N_D(s))\Vert_\infty ds\\ 
& \leq \int_0^t\Vert\nabla_N^+T_N(t-s)\Vert_{L^\infty(I)\rightarrow L^\infty(I)}|F(\bar u^N_C(s),u^N_D(s))|ds \\
& \leq cM_F\int_0^t(t-s)^{-1/2}ds\\ 
& \leq cM_F\sqrt{t} 
\end{align*}
Recall that $\displaystyle b^N(t)\in\mathds{H}^N$, it is a piecewise constant function of the spatial variable. We approximate it by a continuous function of the spatial variable $\displaystyle \tilde{b}^N(t)\in C(I)$, 
using piecewise linear interpolation. Clearly:
\[ \left\Vert \tilde{b}^N(t)\right\Vert_\infty=\left\Vert b^N(t)\right\Vert_\infty. \] 
Since $\tilde{b}^N(t)$ is piecewise linear with $\displaystyle\nabla \tilde{b}^N(t)= \nabla_N^+ b^N(t)$, we deduce form the above estimates that for each 
$t$, $\tilde{b}^N(t)$  is bounded in the space of Lipschitz functions on $I$ and is therefore compact in $C(I)$.

It remains to write 
$$
\left\Vert b^N(t)-\tilde{b}^N(t)\right\Vert_\infty \le \frac1N \left\Vert \nabla_N^+b_C^N(t)\right\Vert_\infty\longrightarrow0, 
$$
as $N\rightarrow\infty$ to get the compactness of $\displaystyle \left\{b^N\right\}_N$ in $D(I)$

\medskip

It finally follows that $\displaystyle \left\{\tilde{u}_C^N\right\}_N$ is $C-$tight and as claimed the limit points are in $C([0,T],C(I))$.

\subsubsection{Tightness for the discrete component}\label{tightness_for_the_discrete_component}

We prove that $\displaystyle \left\{u_D^N\right\}_N$ is tight. Let $(\Omega,\mathcal{F},\mathds{P})$ be the abstract probability space where our stochastic processes $u^N$ are defined. Let $\displaystyle \big\{Y_{j+},Y_{j-},Y_{j,r},Y_{k,r},1\leq j\leq N,1\leq k\leq m,r\in\mathfrak{R}\big\}$ be a family of independent standard Poisson processes. We know that for all $N$, the Markov process $\displaystyle u^N=\big\{u^N(t),t\geq0\big\}$ is progressive, since it is c\`adl\`ag. Moreover, its generator is  $\mathcal{A}^N$. Thus, by Proposition 1.7, Part 4 of \cite{Kurtz1986}, $\displaystyle u^N$ is solution to the martingale problem associated with $\mathcal{A}^N$ (on the specified domain), in the sense that for all  measurable and bounded 
$\displaystyle \varphi$, the process $M_\varphi$ defined by
\begin{equation}\label{martingale_problem_definition_for_the_sequence}
\displaystyle M_\varphi^N(t)=\varphi\big(u^N(t)\big)-\varphi\big(u^N(0)\big)-\int_0^t\mathcal{A}^N\varphi\big(u^N(s)\big)ds
\end{equation} 
for all $t\geq0$, is a $\mathds{P}$-martingale with respect to the (usual) natural associated filtration. 


In addition, $\mathcal{A}^N$ is bounded on the domain  consisting of bounded measurable functions, since reaction rates are bounded. Hence, by Theorem 4.1, Part 6, of \cite{Kurtz1986}, we know that, there exists a sequence $\displaystyle \left(\hat{u}^N=\big(\hat{u}_C^N,\hat{u}_D^N\big)\right)_N$ of c\`adl\`ag stochastic processes in $\displaystyle D\left([0,T],D(I)\times \R^k\right)$  such that they have the same laws as $u^N=(u^N_C,u^N_D)$ and \\

\noindent $\displaystyle {\hat u}_D^N(t) = {\hat u}_D^N(0)+\sum_{\ell=1}^k\left\{\sum_{r\in\mathfrak{R}_{DC}\backslash S_1}
\gamma_{r}^D \1^\ell Y_{\ell,r}\left(\int_0^t\lambda_r\left(\sum_{j=N_{\ell-1}+1}^{N_\ell}a_{j,N}^r{\hat u}_{j}^{N,C}(s),{\hat u}_\ell^{N,D}(s)\right)ds\right)\right.$\par 

$\displaystyle \hspace{3.5cm}\left.+\sum_{r\in\mathfrak{R}_D}\gamma_{r}^D\1^\ell Y_{\ell,r}\left(\int_0^t\lambda_r\left({\hat u}_\ell^{N,D}(s)\right)\right)ds\right\rbrace$.\vspace{0cm}\\
and a similar expression hold for $\hat u_C^N$.
Since we are only interested in the laws of the processes, we omit the hat. 

Let us show that for each $\ell, r$, the laws of $Z_{\ell,r}^N(t)=Y_{\ell,r}(\int_0^t \tilde\lambda_r^N(s)ds)$ for $N\in \N$ are tight in $\displaystyle D\big([0,T],\R\big)$ where 
$\tilde\lambda_r^N(s)$ is either $\lambda_r\left(\sum_{j=N_{\ell-1}+1}^{N_\ell}a_{j,N}^ru_{j}^{N,C}(s),u_\ell^{N,D}(s)\right)$ or $\lambda_r\left(u_\ell^{N,D}(s)\right)$.

From Theorem 1.3, Section 1, Chapter 1 of \cite{Billingsley1999}, we know that the law $Y_{\ell,r}$ is tight in $D([0,T];\R)$ for all $T\ge 0$. Hence, for every $\varepsilon>0$, there exists a 
compact set $K_\varepsilon^T$ in $D([0,T];\R)$ such that 
$$
\P(Y_{\ell,r}\in K_\varepsilon^T)\ge 1-\epsilon.
$$
From Theorem 6.3, Remark 6.4, part 3 in \cite{Kurtz1986}, we may assume that 
$$
K_\varepsilon^T=\{f\in D([0,T];\R): |f(t)|\le M_\varepsilon^T, \; w'(f,\delta,T)\le \alpha_\varepsilon^T(\delta)\mbox{ for all } \delta>0\},
$$ 
with $\lim_{\delta\to0} \alpha_\varepsilon^T(\delta)=0$. The modulus of continuity $w'$ is defined in section 6, part 3 of \cite{Kurtz1986}. Since $w'(Z_{\ell,r}^N,\delta,T)\le w'(Y_{\ell,r},\bar \lambda\delta,\bar\lambda T )$, where $\bar\lambda$ is an upper bound of all reaction rates,
we deduce that 
$$
\P(Z_{\ell,r}^N\in \tilde K_\varepsilon^T)\ge 1-\epsilon
$$
with
$$
\tilde K_\varepsilon^T= \{f\in D([0,T];\R): |f(t)|\le M_\varepsilon^T, \; w'(f,\delta,T)\le \alpha_\varepsilon^{\bar \lambda T}(\bar \lambda\delta)\mbox{ for all } \delta>0\}
$$
which is a compact set in $D([0,T];\R)$. 

Since $u^N_D(0)$ converges in probability, we deduce that the laws of $u^N_D$ for $N\in \N$ are tight in $\displaystyle D\big([0,T],\R^k\big)$.

\subsubsection{Conclusion about tightness} From Corollary 3.33, Section 3, Chapter 6 of \cite{Jacod1987}, $\displaystyle \left\{u^N=\big(u_C^N,u_D^N\big)\right\}_N$ is tight in $\displaystyle D\big([0,T],D(I)\big)\times D\big([0,T],\R^k\big)\equiv D\big([0,T],D(I)\times \R^k\big)$. In other words, $\displaystyle \big\{P_N\big\}_N$ is tight in $\displaystyle \mathcal{P}\big(D([0,T],D(I)\times \R^k)\big)$. 
%
%
%

\subsection{Identification of the limit}

For all $N$, $u^N$ is a Markov jump process, and thus, is solution to the martingale problem associated with its generator $\mathcal{A}^N$, in the sense given by $(\ref{martingale_problem_definition_for_the_sequence})$. In particular,  for all $\displaystyle \varphi\in \cE$ for all $n\geq 1$, for all $0\leq t_1,\cdots,t_n\leq s\leq t\leq T$, and for all $\displaystyle \psi\in C_b((D(I)\times\R^k)^n)$,

$$
\begin{array}{l}
\disp \mathds{E}_\mathds{P}\left[M_\varphi^{N}(t)\psi\big(u^N(t_1),\cdots,u^N(t_n)\big)\right]\vspace{0.2cm}\\

\disp\hspace{2cm} =\mathds{E}_\mathds{P}\left[M_\varphi(s)^N\psi\big(u^N(t_1),\cdots,u^N(t_n)\big)\right],
\end{array}
$$
where\par 
$\disp\hspace{2cm} M_\varphi^N(t)=\varphi\big(u^N(t)\big)-\varphi\big(u^N(0)\big)-\int_0^t\mathcal{A}^N\varphi\big(u^N(s)\big)ds$.\\

\noindent Since the family $\displaystyle \big\{P_N\big\}_N$ is tight in $\displaystyle \mathcal{P}\big(D([0,T],D(I)\times \R^k)\big)$, it is relatively compact there, by Prohorov theorem. Therefore, there exists a subsequence $\displaystyle \big(P_{N_l}\big)_l$ and a probability measure $\displaystyle P\in\mathcal{P}\big(D([0,T],D(I)\times \R^k)\big)$, such that \[ P_{N_l}\Longrightarrow P\hspace{0.3cm} \text{(weakly) as}\hspace{0.3cm} l\rightarrow \infty. \] 
Equivalently, there exists a process $u$ with sample paths in $\displaystyle D\big([0,T],D(I)\times \R^k\big)$ and whose law is $P$, such that the subsequence $\big(u^{N_l}\big)_l$ satisfies $u^{N_l}\Longrightarrow u$ as $l\rightarrow\infty$. The law of $\displaystyle u^{N_l}$ is denoted by $P_{N_l}$ for each $l$. The induced subsequence of martingale problems reads 

\begin{equation}\label{martingale_problem_expectation_representation_for_the_subsequence}
\begin{array}{l}
\disp \mathds{E}_\mathds{P}\left[M_\varphi^{N_l}(t)\psi\big(u^{N_l}(t_1),\cdots,u^{N_l}(t_n)\big)\right]\vspace{0.2cm}\\

\disp\hspace{2cm}=\mathds{E}_\mathds{P}\left[M_\varphi^{N_l}(s)\psi\big(u^{N_l}(t_1),\cdots,u^{N_l}(t_n)\big)\right].
\end{array}
\end{equation}\vspace{0cm}

There are several difficulties if we try to take the limit in \eqref{martingale_problem_expectation_representation_for_the_subsequence}, as $l\rightarrow\infty$. 
In particular, $M_\varphi^{N_l}$ contains the diffusion term which is linear in $u^N$. This creates difficulties when 
taking limits inside the expectation since we do not have any estimates on the moments of $u^N$. Tightness was obtained through bounds in probability.

To avoid this problem, we consider the process $w^N=(v^N_C,u^N_C,u^N_D)$, where $v^N_C$ was introduced in \cref{s3.1.1}. It is a Markov process. In fact it is a finite dimensional PDMP 
with generator 
\[
\begin{array}{l}
\displaystyle {\widetilde{\mathcal{A}}}^N\varphi(v^N_C,u^N_C, u^N_D)\vspace{0.2cm}\\  
\displaystyle \hspace{0.2cm}=\left\langle D_{v^N_C}\varphi(v^N_C,u^N_C, u^N_D),\Delta_N v^N_C+F(v^N_C,u^N_D)\right\rangle_2\vspace{0.1cm}\\

\displaystyle \hspace{1cm}+\sum_{j=1}^N\left\lbrace\sum_{r\in\mathfrak{R}_C}\left[\varphi\left(v^N_C,u^N_C+
\frac{\gamma_r^C}{\mu}\mathds{1}_j,u^N_D\right)-\varphi(v^N_C,u^N_C,u^N_D)\right]\mu\lambda_r\left(u_j^{N,C}\right)\right.\\ 

\displaystyle \hspace{2.5cm}\left.+\sum_{r\in S_1} \left[\varphi\left(v^N_C,u^N_C+\frac{\gamma_r^C}{\mu}\mathds{1}_j,u^N_D\right)-\varphi(v^N_C,u^N_C,u^N_D)\right]\mu\lambda_r\left(u_j^{N,C},u_{\ell_j}^{N,D}\right)\right\}\\ 

\displaystyle \hspace{1cm} +\sum_{j=1}^N\left\lbrace\left[\varphi\left(v^N_C,u^N_C+\frac{\mathds{1}_{j-1}-\mathds{1}_j}{\mu},u^N_D\right) 
- \varphi(v^N_C,u^N_C,u^N_D)\right]\mu N^2u_j^{N,C}\right.\\

\displaystyle \hspace{3cm} +\left.\left[\varphi\left(v^N_C,u^N_C+\frac{\mathds{1}_{j+1}-\mathds{1}_j}{\mu},u^N_D\right)-2\varphi(v^N_C,u^N_C,u^N_D)\right]\mu N^2u_j^{N,C}\right\}\\


\displaystyle \hspace{1cm}+\sum_{\ell=1}^k\left\{\sum_{r\in\mathfrak{R}_{DC}\backslash S_1}\left[\varphi\left(v^N_C,u^N_C+\frac{\gamma_{r}^{C,N}}{\mu},u^N_D+\gamma_{r}^{D}\1^\ell\right)-\varphi(v^N_C,u^N_C,u^N_D)\right]\right.\\
\disp \hfill  \times \lambda_r\left(\sum_{j=N_{\ell-1}+1}^{N_\ell}a_{j,N}u_j^{N,C},u^{N,D}_\ell\right) \\

\displaystyle \hspace{2.5cm}\left.+\sum_{r\in\mathfrak{R}_D}\left[\varphi\left(v^N_C,u^N_C,u^N_D+\gamma_{r}^{D}\1^\ell\right)-\varphi(v^N_C,u^N_C,u^N_D)\right]\lambda_r\left(u_\ell^{N,D}\right)\right\rbrace.
\end{array}
\]

The laws of $w^N$ are tight in $D([0,T];D(I)\times D(I)\times \R^k)$. By the discussion of the begining of \cref{s3.1.1}, we know that the limit points are concentrated on the set 
$\{(u_C,v_C,u_D) \in C([0,T];C(I)\times C(I)) \times D([0,T];\R^k):\? u_C=v_C\}$. We repeat the argument above for $w^N$ and deduce the existence of a subsequence such that  $w^{N_l}\Longrightarrow w$ as $l\rightarrow\infty$

Using the representation theorem of Skorohod, there exist versions of the stochastic processes, $\tilde{w}^{N_l}$, $\tilde{w}$ on a probability space $\big(\tilde{\Omega},\tilde{\mathcal{F}},\tilde{\mathds{P}}\big)$, such that \[ \tilde{w}^{N_l}\longrightarrow \tilde{w} \hspace{0.5cm}\tilde{\mathds{P}}-a.s.\footnote{abbreviation for almost surely.} \] in the Skorohod topology. Also, we know that $\tilde w$ is of the form $(\tilde u_C,\tilde u_C,\tilde u_D)$ for some processes $\tilde u_C$, $\tilde u_D$ in 
$C([0,T];C(I))$ and $D([0,T];\R^k)$. Below, we write $\tilde v=(\tilde u_C,\tilde u_D)$.

Since we are interested in the laws of the processes, we consider these new versions in the sequel, and conserve the initial notations (without the "tilde"). Moreover, there exists a subset $\underline{\Omega}\subset\Omega$ such that $\mathds{P}(\underline{\Omega})=1$ and for all $\omega\in \underline{\Omega}$, $w^{N_k}(\omega)\longrightarrow w(\omega)$ in the Skorohod topology. \par 

In addition, we know that $\Vert w^{N_l}(t)-w(t)\Vert_\infty\rightarrow 0$  if $w$ is continuous at $t\in [0,T]$ (see Section 12, Chapter 3 of \cite{Billingsley1999}). Thus this holds almost everywhere in $(t,\omega)\in [0,T]\times \Omega$. Also, from Lemma 1, Section 12, Chapter 3 of \cite{Billingsley1999}, $w$ is $\mathds{P}$-a.s. continuous at every $t$, except for a countable set $D_P=D_\mathds{P}(w)$. The set $D_P$ is the complementary of \[ T_P=T_\mathds{P}(w):=\left\{t\in [0,T]:\mathds{P}\big(w(t)= w(t^-)\big)=1\right\}. \] Thus, for all $t\in T_P$, there is a subset $\Omega_t\subset\Omega$ with $\mathds{P}(\Omega_t)=1$, such that $w(\omega)$ is continuous at $t$ for all $\omega\in\Omega_t$. Note that $0,T\in T_P$.\par 

To consider the limit $l\to \infty$, we consider particular test functions depending only on $v^N=(v^N_C,u^N_D)$. For such functions we have the following expression 
for the generator:
\begin{equation}\label{e3.5}
\begin{array}{l}
\displaystyle {\widetilde{\mathcal{A}}}^N\varphi(v^N_C,u^N_C, u^N_D)\vspace{0.2cm}\\  
\displaystyle \hspace{0.2cm}=\left\langle D_{v^N_C}\varphi(v^N_C, u^N_D),\Delta_N v^N_C+F(v^N_C,u^N_D)\right\rangle_2\vspace{0.1cm}\\

\displaystyle \hspace{0.5cm} +\sum_{\ell=1}^k\left\{\sum_{r\in\mathfrak{R}_{DC}\backslash S_1}\left[\varphi\left(v^N_C,u^N_D+\gamma_{r}^{D}\1^\ell\right)-\varphi(v^N_C,u^N_D)\right]\lambda_r\left(\sum_{j=N_{\ell-1}+1}^{N_\ell}a_{j,N}^ru_j^{N,C},u^{N,D}_\ell\right)\right. \\

\displaystyle \hspace{2cm}\left.+\sum_{r\in\mathfrak{R}_D}\left[\varphi\left(v^N_C,u^N_D+\gamma_{r}^{D}\1^\ell\right)-\varphi(v^N_C,u^N_D)\right]\lambda_r\left(u_\ell^{N,D}\right)\right\rbrace
\end{array}
\end{equation}


In particular,  for all $\displaystyle \varphi\in \cE$ for all $n\geq 1$, for all $0\leq t_1,\cdots,t_n\leq s\leq t\leq T$, and for all $\displaystyle \psi\in C_b((D(I)\times\R^k)^n)$,

\begin{equation}\label{e3.6}
\begin{array}{l}
\disp \mathds{E}_\mathds{P}\left[\widetilde M_\varphi^{N}(t)\psi\big(v^N(t_1),\cdots,v^N(t_n)\big)\right]\ =\mathds{E}_\mathds{P}\left[\widetilde M_\varphi(s)^N\psi\big(v^N(t_1),\cdots,v^N(t_n)\big)\right],
\end{array}
\end{equation}
where\par 
$\disp\hspace{2cm}\widetilde M_\varphi^N(t)=\varphi\big(v^N(t)\big)-\varphi\big(v^N(0)\big)-\int_0^t\widetilde{\mathcal{A}}^N\varphi\big(w^N(s)\big)ds$.\\
This follows from the fact that functions in $\cE$  are in the domain of the finite dimensional PDMP $w^N$ (see \cite{Arnaud2012}).

For $t_1,\cdots,t_n,s,t\in T_P$, the set 
\[ \displaystyle \widetilde\Omega=\left(\bigcup_{i=1}^n\Omega_{t_i}\right)\cup\big(\Omega_s\cup\Omega_t\cup \underline{\Omega} \big) \] is of probability $1$, and,  for $\omega\in \tilde\Omega$, we have for $l\rightarrow\infty$:
\begin{itemize}
\item $\displaystyle \psi\big(v^{N_l}(t_1),\cdots,v^{N_l}(t_n)\big)\to \psi\big(v(t_1),\cdots,v(t_n)\big)$, since $\displaystyle v^{N_k}(t_i)\rightarrow_l v(t_i)$ for $i=1,\cdots,n$.

\item $\displaystyle v^{N_l}(0)\to v(0)$ and $\displaystyle \varphi\big(v^{N_l}(0)\big)\to \varphi\big(v(0)\big)$, and, by  (ii) of \cref{hybrid_simplification}, $v(0)=v_0$.

\item $\displaystyle v^{N_l}(t)\to v(t)$ and $\displaystyle \varphi\big(v^{N_l}(t)\big)\to \varphi\big(v(t)\big)$. The same hold when $t$ is replaced by $s$.
\end{itemize}

We cannot let $l\to\infty$ directly in \eqref{e3.6}. Indeed $\widetilde{ \mathcal{A}}^{N_l}\varphi\big(w^{N_l}(t)\big)$ contains $\Delta_{N_l}v_C^{N_l}(t)$
which converges to $\Delta v_C(t)$ but in bad topologies since the convergence of $v^{N_l}_C$ holds only in $D(I)$.  
In order to overcome this difficulty, we first use regularized test functions.\\ 

\noindent \textbf{\textit{Regularization.}} For all $\varphi\in\cE$ and $\varepsilon>0$, we introduce the linear operator $A_{\varepsilon}$ on $D(I)$, and the function $\displaystyle \varphi_{\varepsilon} : D(I)\times \R^k\longrightarrow\mathds{R}$, defined by 
\begin{equation}\label{smoothing_linear_operator_and_test_function}
\disp A_{\varepsilon}:=\big(I_d-\varepsilon\Delta\big)^{-2}\hspace{0.5cm}\text{and}\hspace{0.5cm}\disp \varphi_{\varepsilon}(u_C,u_D):=\varphi(A_{\varepsilon}u_C,u_D),
\end{equation}
%
%
%
It is well known that $A_{\varepsilon}$ is a bounded linear operator of contraction (see e.g. \cite{Cazenave1998}, Chapter 2). Also, it commutes with $\Delta$ and with $\Delta_N$. Moreover, it maps $L^\infty(I)$, and in particular $D(I)$, into $W^{4,\infty}(I)$.

\begin{lemma}\label{properties_of_the_perturbated_test_functions} \hspace{1cm}\vspace{0.1cm}\par 
 For all $g_N,g\in D(I)$ such that $\Vert g_N-g\Vert_\infty\longrightarrow 0$ as $N\rightarrow\infty$,\vspace{0.1cm}\par 
\hspace{1cm} (i) $\disp \left\Vert A_{\varepsilon}g_N-A_{\varepsilon}g\right\Vert_{C^2(I)}\longrightarrow0$ as $N\rightarrow\infty$, \vspace{0.1cm}\par 
\hspace{1cm} (ii) $\disp \Vert \Delta_NA_{\varepsilon}g_N-\Delta A_{\varepsilon}g\Vert_\infty\longrightarrow 0$ as $ N\rightarrow\infty$.\vspace{0.2cm}\par 
\end{lemma}

The second point (ii) immediately follows from the first (i), thanks to $\cref{Discrete_Laplacian_properties}$ (vii). Point  (i) 
is a consequence of the fact that $A_{\varepsilon}$ is a continuous operator from $L^\infty(I)$ to $W^{4,\infty}$ and therefore to 
$C^2(I)$.

Let us 
consider $(\ref{e3.6})$, with the test function $\varphi_{\varepsilon}$ instead of $\varphi$:
\begin{equation}\label{pertubated_martingale_problem_expectation_representation_for_the_subsequence}
\begin{array}{l}
\disp \mathds{E}_\mathds{P}\left[\widetilde M_{\varphi_{\varepsilon}}^{N_l}(t)\psi\big(v^{N_l}(t_1),\cdots,v^{N_l}(t_n)\big)\right]\vspace{0.2cm}=\mathds{E}_\mathds{P}\left[\widetilde M_{\varphi_{\varepsilon}}^{N_l}(s)\psi\big(v^{N_l}(t_1),\cdots,v^{N_l}(t_n)\big)\right],
\end{array}
\end{equation}
where\par 
$\disp \hspace{1.5cm} \widetilde M_{\varphi_{\varepsilon}^{N_l}}(t)=\varphi_{\varepsilon}\big(v^{N_l}(t)\big)-\varphi_{\varepsilon}\big(v^{N_l}(0)\big)-\int_0^t\widetilde{\mathcal{A}}^{N_l}\varphi_{\varepsilon}\big(w^{N_l}(r)\big)dr$.\\
 
\noindent We start keeping $\varepsilon$ fixed, and let $l\rightarrow\infty$. Clearly, $\varphi_{\varepsilon}$ is continuous on $D(I)$. Since $s,t\in T_P$,  it follows from the preceding discussion that 
\[ \left\{
\begin{array}{l}
\displaystyle \varphi_{\varepsilon}\big(v^{N_l}(0)\big)\longrightarrow_l \varphi_{\varepsilon}(v_0)\hspace{0.2cm}a.s.\vspace{0.1cm}\\

\displaystyle \varphi_{\varepsilon}\big(v^{N_l}(s)\big)\longrightarrow_l \varphi_{\varepsilon}\big(v(s)\big)\hspace{0.2cm}a.s.\vspace{0.1cm}\\

\displaystyle \varphi_{\varepsilon}\big(v^{N_l}(t)\big)\longrightarrow_l \varphi_{\varepsilon}\big(v(t)\big)\hspace{0.2cm}a.s.
\end{array} 
\right. \] 
It is rather straightforward to prove that $\disp\widetilde{\mathcal{A}}^{N_l}\varphi_{\varepsilon}\big(w^{N_l}(t)\big)\longrightarrow_l \mathcal{A}^\infty\varphi_{\varepsilon}\big(v(t)\big)$. Also, it is uniformly bounded in $N$ on $[0,T]$. This follows form the boundedness of the reaction rates, the bound \eqref{e3.2} and the fact that
$A_\varepsilon$ is a bounded operator on $D(I)$.

Therefore, by dominated convergence, we may let $l\to \infty$ in \cref{pertubated_martingale_problem_expectation_representation_for_the_subsequence} and obtain
\begin{equation}\label{pertubated_martingale_problem_expectation_representation}
\begin{array}{l}
\disp \mathds{E}_\mathds{P}\left[M_{\varphi_{\varepsilon}}^{\infty}(t)\psi\big(v(t_1),\cdots,v(t_n)\big)\right] =\mathds{E}_\mathds{P}\left[M_{\varphi_{\varepsilon}}^{\infty}(s)\psi\big(v(t_1),\cdots,v(t_n)\big)\right],
\end{array}
\end{equation}
with 
\begin{equation}\label{limiting_regularized_martingale}
\displaystyle M_{\varphi_{\varepsilon}}^\infty(t)=\varphi_{\varepsilon}\big(v(t)\big)-\varphi_{\varepsilon}(v_0)-\int_0^t\mathcal{A}^\infty\varphi_{\varepsilon}\big(v(r)\big)dr,
\end{equation}
where $\displaystyle \mathcal{A}^\infty$ is the generator defined by $(\ref{limiting_generator})$.

We now want to let $\varepsilon\to 0$. Except for the one containing the Laplace operator, all terms in the \eqref{pertubated_martingale_problem_expectation_representation} 
are easily seen to converge.   To treat the remaining term, we observe that we may take the limit in \cref{e3.2bis} along the subsequence $N_l$ and 
deduce that $u_C$ satisfies almost surely:
$$
u_C(t)=T(t)v_{C,0}+\int_0^t T(t-s)F(u_C(s),u_D(s))ds,\, a.s.
$$
The easiest way to do this is to take the limit in the weak form of the equation. From the smoothing property of $T(t)$:
$$
\|T(t)\|_{L^\infty(I)\to W^{1,\infty}(I)}\le c (t^{-\frac12}+1),\; t\ge 0, 
$$
we deduce that $u_C(t)$ is bounded uniformly in $W^{1,\infty}(I)$, and in particular in $H^{1}_0(I)$ for $t\in [0,T]$ and $\omega \in \underline{\Omega}$. We deduce
that $\Delta u_C(t)$ is bounded uniformly in $H^{-1}$ and thanks to the property of the differential of functions in $\cE$ we may let $\varepsilon\to 0$ in the
term containing the Laplace operator.

We obtain for $\varphi\in\cE$:
\begin{equation}\label{martingale_problem_expectation_representation}
\begin{array}{l}
\disp \mathds{E}_\mathds{P}\left[M_{\varphi}^{\infty}(t)\psi\big(v(t_1),\cdots,v(t_n)\big)\right]=\mathds{E}_\mathds{P}\left[M_{\varphi}^{\infty}(s)\psi\big(v(t_1),\cdots,v(t_n)\big)\right],
\end{array}
\end{equation}
with
\begin{equation}\label{limiting_martingale}
\displaystyle M_{\varphi}^\infty(t)=\varphi\big(v(t)\big)-\varphi(v_0)-\int_0^t\mathcal{A}^\infty\varphi\big(v(r)\big)dr.
\end{equation} 

%
%

Now, if any of $t_1,\cdots,t_n,s,t$ does not lie in $T_P$, let us say $t\notin T_P$ for instance, we choose a sequence $\displaystyle \big(t^l\big)_l$ in $T_P$ such that $t^l\rightarrow_l t$ with $t^l>t$. Since $v$ is \textit{c\`adl\`ag}, it is right-continuous at $t$, and $\displaystyle v\big(t^l\big)\longrightarrow_l v(t)$. Then, we use $(\ref{limiting_martingale})$ with $t^l$ instead of $t$, let $l\rightarrow\infty$ and deduce that $(\ref{limiting_martingale})$ also holds for $t\notin T_P$.

We have proved that the probability measure $P$, the law of $v$, is a solution of the martingale problem associated with the generator $\disp \mathcal{A}^\infty$ on the domain $\cE$. 
Since the reaction rates are bounded as well as their derivatives,  \cref{well_posedness_of_the_martingale_problem_for_inf_dim_PDMP} holds and the martingale problem for $\cA^\infty$ ---restricted to $\cE$--- admits a unique solution, which is the law $P_{v_0}$ of the PDMP $v$ characterized by $(\Delta,F,\lambda^D,Q)$, and which starts at $v_0$. It follows that $P=P_{v_0}$, and the whole sequence of the laws of $v^N$   converges to $P_{v_0}$. This implies the convergence of $(P_N)_N$ to $P_{v_0}$.

\subsection{Conclusion}
\label{s3.3}

Now, we get rid of the additional assumption of boundedness of the process, 
and prove \cref{hybrid_simplification} by a truncation argument. \par 

Let $\eta\in C^{\infty}(\R_+)$ such that
\[ \left\{ 
\begin{array}{l}
\eta(y)=1,\hspace{0.5cm}y\in [0,1],\vspace{0.1cm}\\
\eta(y)=0,\hspace{0.5cm}y\in [2,\infty).
\end{array}
\right. \]
For $n\geq1$ and $r\in\fR$, define  \[ \eta_n(y)=\eta\left(\frac{|y|^2}{n^2}\right)\hspace{0.5cm}\text{and}\hspace{0.5cm}\lambda_r^n(y)=\eta_n(y)\lambda_r(y) \] for $y\in \R^2$. Since (we have in particular) $y\mapsto|y|^2$ belongs to $C^1(\R^2)$ and $\eta\in C_b^1(\R_+)$, we also have $\eta_n\in C_b^1(\R_+)$. Furthermore, $\eta_n$ and its derivative vanish outsite the compact $\bar{B}(0,n)$, which is the closed ball in $\R^2$ centered at $0$ and of radius $n$. Then, the problem with $\lambda_r^n$ instead of $\lambda_r$ fulfills the additional asumptions of the previous steps.\par 

We define $u_n^N=\big(u_{C}^{N,n},u_{D}^{N,n}\big)$ the (truncated) jump Markov process associated to the (truncated) jump intensities $\lambda_r^n$, and starting at $u^N(0)$. 
By the preceding results, we know that, for all $n\geq1$, \[ u^{N,n}\Longrightarrow_{N} v^n,\hspace{0.2cm}\text{in }D(\R_+;E), \] where $v^n$ is the (truncated) PDMP whose (truncated) characteristics $(\Delta,F_n,\lambda_n^D,Q_n)$ are obviously defined, w.r.t. the truncation. 

It remains to argue as in \cite{Arnaud2012} {in order to conclude and to end the proof} of \cref{hybrid_simplification}. $\blacksquare$

\section{Well-posedness of the martingale problem for the generator of an infinite dimensional PDMP}\label{proof_of_well_posedness_of_the_martingale_problem_for_inf_dim_PDMP}
\label{s4}

\noindent \textbf{\underline{Proof of \cref{well_posedness_of_the_martingale_problem_for_inf_dim_PDMP}}.} 
Let $(P_t):=(P_t)_{t\geq0}$ be the semigroup on $E=C(I)\times \R^k$ associated to the PDMP $\disp \left\{v(t)=\big(v_C(t),v_D(t)\big),t\geq0\right\}$ starting at $(v_{C,0},v_{D,0})$
We need the following:

\begin{lemma}\label{stability_of_test_function_space_by_the_semigroup}
For all $t\geq0$, $\varphi\in\cE$, $P_t\varphi$ is bounded differentiable  with respect to the first varibale on $C(I)$ and satisfies: 
\begin{align}
\disp \Vert P_t\varphi\Vert_\infty\leq \Vert \varphi\Vert_\infty,\hspace{6.65cm}\label{contraction_property}\vspace{0.1cm}\\
\disp |D_\alpha P_t\varphi(\alpha,\nu)\cdot h|\leq c\Vert\varphi\Vert_{\cE}e^{Mt}\Vert h\Vert_\infty, \hspace{0.5cm}\alpha\in C(I),\nu\in \R^k, h\in C(I)\label{Lipschitz_property_and_constant}
\end{align}
\noindent for some constants $c$ and $M$ depending only on the characteristics of the PDMP $v$.
\end{lemma}

\noindent \textit{Proof of $\cref{stability_of_test_function_space_by_the_semigroup}$}. The part $(\ref{contraction_property})$ immediately follows, from the definition of the semigroup and the fact that the expectation is increasing. For $\varphi\in\cE$ and $\psi\in\fB(\R_+\times E)$, let us define
\begin{align*}
\disp G_\varphi\psi(t,u)&:=\E_{u}\left[\varphi(v(t))\1_{t<T_1}+\psi(t-T_1,v(T_1))\1_{t\geq T_1}\right]\\
									  &= \varphi(\phi_\nu(t,\alpha),\nu)H(t,u)\\
									  & \hspace{0.5cm}+\int_{0}^t\int_K\psi(t-s,\alpha,\xi)Q(d\xi;\phi_\nu(s,\alpha),\nu)\Lambda(\phi_\nu(s,\alpha),\nu)H(s,u)ds
\end{align*}
for $(t,u)\in\R_+\times E$, with $u=(\alpha,\nu)$ and $v$ the PDMP starting at $u$. Then, according to Lemma 27.3 of \cite{Davis1993}, 

\[ G_\varphi^n\psi(t,u)=\E_u\left[\varphi(v(t))\1_{t<T_n}+\psi(t-T_n,v(T_n))\1_{t\geq T_n}\right] \] and \[ \lim_{n\rightarrow\infty}G_\varphi^n\psi(t,u)=P_t\varphi(u), \] 

\noindent where
\begin{align*}
\disp G_\varphi^n\psi(t,u) & :=G_\varphi\big(G_\varphi^{n-1}\psi(t,u)\big)\\
\disp 									& = \varphi(\phi_\nu(t,\alpha),\nu)H(t,u)\\
\disp 									& \hspace{0.5cm}+\int_{0}^tH(s,u)\left[\int_K\Lambda_\nu(\phi_\nu(s,\alpha))G_\varphi^{n-1}\psi(t-s,\alpha,\xi)Q(d\xi;\phi_\nu(s,\alpha),\nu)\right]ds.
\end{align*}

Thus, our assumptions allow us to use dominated convergence, and we deduce
\begin{equation}\label{semigroup_as_fix_point_of_the_functional_G}
\begin{array}{l}
\disp P_t\varphi(u)=\varphi(\phi_\nu(t,\alpha),\nu)H(t,u)\vspace{0.1cm}\\
\disp \hspace{2cm}+\int_{0}^t\int_KP_{t-s}\varphi(\alpha,\xi)Q(d\xi;\phi_\nu(s,\alpha),\nu)\Lambda(\phi_\nu(s,\alpha),\nu)H(s,u)ds.
\end{array}
\end{equation}
\noindent That is, $P_\cdot\varphi:t\mapsto P_t\varphi$ is a fixed point of $G_\varphi$.\par 

For $T>0$, we introduce the Banach space $\cE_T:=L^\infty([0,T],C_b^{1,0}(C(I)\times \R^k)$, with the norm \[ \Vert \psi\Vert_{\cE_T}=\sup_{t\in [0,T]}e^{-\beta t}(\Vert\psi(t,\cdot)\Vert_{\infty}+\|D_\alpha \psi(t,\cdot)\|_\infty), \] where $\beta$ will be fixed hereafter.

For $t\in[0,T]$, $(\alpha,\nu)\in E$, set \vspace{0.1cm}\par 
$\disp\hspace{5cm} q_\varphi(t,\alpha,\nu)=\varphi(\phi_\nu(t,\alpha),\nu)H(t,\alpha,\nu)$.\\ We claim that 
\begin{equation}\label{intermediate_functional_q}
\disp q_\varphi:t\mapsto q_\varphi(t,\cdot)\hspace{0.2cm} \text{is in}\hspace{0.2cm} \cE_T.
\end{equation}\vspace{0.1cm}

\noindent Indeed, let $0\leq t\leq T$, $u=(\alpha,\nu)\in E$ be fixed. Then \[ |H(t,\alpha,\nu)|\leq 1, \] from the definition of the survivor function $H$ by $(\ref{survivor_function_of_the_first_transition_time})$. Since $\varphi$ is bounded, we have \[ \Vert q_\varphi(t,\cdot)\Vert_\infty \leq \Vert\varphi\Vert_\infty\hspace{0.5cm}\text{and}\hspace{0.5cm} q_\varphi(t,\cdot)\in\fB_b(E). \] Next, we show that $q_\varphi(t,\cdot)\in\cE$. Let $D_\alpha$ denote the differential operator w.r.t. the variable $\alpha$, and set \[ \Psi_\nu(t)\cdot h=D_\alpha\phi_\nu(t,\alpha)\cdot h. \] 

\noindent From our assumptions, $F\in C_b^{1,0}(E)$. It follows that for $h\in C(I)$ the map $t\mapsto\Psi_\nu(t)\cdot h$ is a global mild solution of
\begin{equation}\label{differential_of_the_flow_wrt_the_initial_condition_equations}
\left\{
\begin{array}{l}
\disp \frac{\partial}{\partial t}\Psi_\nu(t)\cdot h=\Delta\Psi_\nu(t)\cdot h+D_\alpha F_\nu\big(\phi_\nu(t,\alpha)\big)\circ\Psi_\nu(t)\cdot h \vspace{0.1cm}\\
\disp \Psi_\nu(0)\cdot h=h.
\end{array} 
\right.
\end{equation}
It satisfies \[ \Psi_\nu(t)\cdot h=T(t)h+\int_0^t T(t-s)D_\alpha F_\nu\big(\phi_\nu(s,\alpha)\big)\circ\Psi_\nu(s)\cdot h ds. \] Since the semigroup $\{T(t)\}$ of $\Delta$ is of contraction,  the Gronwall Lemma yield, \[ \Vert\Psi_\nu(t)\cdot h\Vert_\infty \leq \|h\|_\infty+L_F\int_0^t\Vert\Psi_\nu(s)\cdot h\Vert_\infty ds \leq e^{L_Ft}\|h\|_\infty. \] Thus, \[ \Vert D_\alpha\phi_\nu(t,\alpha)\cdot h\Vert_\infty \leq e^{L_Ft}\|h\|_\infty. \]

The chain rule leads to
\begin{align*}
\disp \Vert D_\alpha H(t,\alpha,\nu)\cdot h\Vert_\infty & \leq \int_0^t\left\Vert D_\alpha\left[\Lambda\big(\phi_\nu(s,\alpha),\nu\big)\right]\cdot h\right\Vert_\infty ds \leq L_\Lambda\int_0^t\Vert D_\alpha\phi_\nu(s,\alpha)\cdot h\Vert_\infty ds\\
\disp & \leq \frac{L_\Lambda}{L_F}\left(e^{L_Ft}-1\right)\|h\|_\infty.
\end{align*}
Therefore $q_\varphi(t,\cdot)\in \cE$ , since
\begin{align*}
\disp \Vert D_\alpha q_\varphi(t,\alpha,\nu)\cdot h \Vert_\infty & \leq \Vert D_\alpha\varphi(\phi_\nu(t,\alpha),\nu)\cdot h \Vert_\infty+\Vert D_\alpha H(t,\alpha,\nu)\cdot h \Vert_\infty\\
\disp & \leq \left(\|D_\alpha \varphi\|_{\infty}+\frac{L_\Lambda}{L_F}\left(e^{L_Ft}-1\right)\right)\|h\|_\infty.
\end{align*} 
Finally, for $\beta>L_F$, 
\begin{equation}\label{bound_for_q_varphi}
\disp \Vert q_\varphi\Vert_{\cE_T} \leq \Vert \varphi\Vert_\infty+\left(\|D_\alpha \varphi\|_{\infty}+\Vert\varphi\Vert_\infty\frac{L_\Lambda}{L_F}\right).
\end{equation}

Similar arguments will allow us to get the upcoming upper bounds. If $\psi\in \cE_T$, then $G_\varphi\psi\in\cE_T$. In fact, \[ G_\varphi\psi(t,\alpha,\nu)=q_\varphi(t,\alpha,\nu)+\tilde{q}_\psi(t,\alpha,\nu), \] with\par  
$\disp\hspace{0.2cm} \tilde{q}_\psi(t,\alpha,\nu)=\int_0^tH(s,\alpha,\nu)\left[\Lambda\big(\phi_\nu(s,\alpha),\nu\big)\int_K\psi(t-s,\alpha,\xi)Q_\nu\big(d\xi;\phi_\nu(s,\alpha)\big)\right]ds$. \\
It is standard to prove that this defines a function in $\cE_T$. 
Moreover, observing that  $\Vert\tilde{q}_\psi(t,\cdot)\Vert_\infty\leq M_\Lambda\Vert\psi\Vert_\infty t\leq c_1t$, and recalling the definition of $\Lambda$ and  $Q$ 
(see \eqref{global_jump_rate_slow_reaction} and \eqref{transition_measure}), which in particular imply their differentiability, 
\begin{align*}
 \Vert D_\alpha\tilde{q}_\psi(t,\alpha,\nu)\cdot h \Vert_\infty & \leq \int_0^t\left[\frac{L_\Lambda}{L_F}\big(e^{L_Fs}-1\big)M_\Lambda\Vert\psi(t-s,\cdot)\Vert_\infty+L_Q\Vert\psi\Vert_{\cE_T}\right]ds \|h\|_\infty,
\end{align*} 
it follows from $(\ref{bound_for_q_varphi})$ and $(\ref{contraction_property})$ that
\begin{align*}
\Vert D_\alpha G_\varphi\psi(t,\alpha,\nu)\cdot h \Vert_\infty& \leq \bigg( \|D_\alpha \varphi\|_{\infty}+\Vert\varphi\Vert_\infty\frac{L_\Lambda}{L_F}e^{L_Ft}\\
&\hspace{0.5cm} +\left(\frac{L_\Lambda}{L_F}M_\Lambda+\frac{L_Q}{L_F}\right)\Vert\varphi\Vert_\infty e^{L_Ft}+L_Q\int_0^tL_{\psi(t-s,\cdot)}ds\bigg) \|h\|_\infty,
\end{align*}
which yields (using Gronwall lemma)
\begin{equation}\label{Lipschitz_constant_determining}
\disp \|D_\alpha G_\varphi \psi(t,\cdot)\|_{\infty}\leq c\Vert\varphi\Vert_{C^{1,0}_b(C(I)\times \R^k)} e^{L_Ft}+L_Q\int_0^t \|D_\alpha\psi(s,\cdot)\|_{\infty} ds\leq c\Vert\varphi\Vert_{C^{1,0}_b(C(I)\times \R^k)} e^{Mt},
\end{equation}
and we conclude that $G_\varphi$ maps $\cE_T$ into itself. The constants $c$ and $M$ depend only on ($L,F,\Lambda,Q$) and will turn out to be the constants appearing in $(\ref{Lipschitz_property_and_constant})$.\vspace{0.1cm}\par 

Moreover, if $\psi_1,\psi_2\in\cE_T$, we prove similarly 
\begin{align*}
\|G_\varphi\psi_1(t,\cdot)-G_\varphi\psi_2(t,\cdot)\Vert_\infty\leq M_\Lambda\int_0^t\Vert\psi_1(t-s,\cdot)-\psi_2(t-s,\cdot)\Vert_\infty ds,\hspace{0.5cm}\\
\Vert D_\alpha G_\varphi\psi_1(t,\cdot)-D_\alpha G_\varphi\psi_2(t,\cdot)\Vert_\infty\leq \kappa_1\int_0^te^{L_Fs}\Vert \psi_1(t-s,\cdot)-\psi_2(t-s,\cdot)\Vert_{C^{1,0}_b(C(I)\times \R^k)} ds,
\end{align*}
where $\kappa_1$ is a constant depending only on ($L,F,\Lambda,Q$).\par 
Then, it is not difficult to deduce that
\[
\begin{array}{l}
\Vert G_\varphi\psi_1-G_\varphi\psi_2\Vert_{\cE_T}\vspace{0.1cm}\\
\hspace{1cm} \leq \kappa_2 \sup_{t\in[0,T]}e^{-\beta t}\left(\int_0^t e^{\beta(t-s)}ds+\int_0^t e^{\beta(t-s)} e^{L_Fs}ds\right)\Vert\psi_1-\psi_2\Vert_{\cE_T}\vspace{0.1cm}\\
\hspace{1cm}\leq \kappa_2\left(\frac{1}{\beta}+\frac{1}{\beta-L_F}\right)\Vert\psi_1-\psi_2\Vert_{\cE_T},
\end{array}
\]
where again, $\kappa_2$ depends only on ($L,F,\Lambda,Q$). We now choose $\beta$ sufficiently large and deduce from the Picard theorem that $G_\varphi$ has a unique fixed point in $\cE_T$. This fixed point is the limit of $G_\varphi^n\psi$ for any $\psi\in\cE_T$. Thus, $P_\cdot\varphi\mapsto P_t\varphi$ is that fixed point.

The Lipschitz constant $L_{P_t\varphi}$ of $P_t\varphi$ is obtained from $(\ref{Lipschitz_constant_determining})$, by taking $\psi=P_\cdot\varphi$. $\square$

\begin{lemma}\label{l4.2}
For all $t\geq0$, $\varphi\in\cE$, $P_t\varphi\in \cE$  and satisfies: 
\begin{align}
\disp |D_\alpha P_t\varphi(\alpha,\nu)\cdot h|\leq \tilde c(t^{-\frac34}+1) e^{\tilde  Mt}\Vert\varphi\Vert_{\cE}  \Vert h\Vert_{H^{-1}}, \hspace{0.5cm}\alpha\in C(I),\nu\in \R^k, h\in C(I)\label{e4.8}
\end{align}
\noindent for some constants $\tilde  c$ and $\tilde  M$ depending only on the characteristics of the PDMP $v$.
\end{lemma}
\noindent \textit{Proof.} If we prove \eqref{e4.8}, it follows that $D_\alpha P_t \varphi$ can be extended to $H^{-1}(I)$ and using \cref{stability_of_test_function_space_by_the_semigroup} the result follows. The proof of \eqref{e4.8} relies on the smoothing properties of the heat kernel. 
Indeed, we have for $t>0$
$$
\|T(t)\|_{H^{-1}(I)\to L^\infty(I)} \le c(t^{-\frac34}+1)
$$
for some constant $c>0$. It follows:
\[ 
\|\Psi_\nu(t)\cdot h\|_\infty\le c(t^{-\frac34}+1)\|h\|_{H^{-1}(I)} + L_F \int_0^t  \| \Psi_\nu(s)\cdot h\|_\infty ds \] 
and by Gronwall Lemma:
\[ \|\Psi_\nu(t)\cdot h\|_\infty\le c(t^{-\frac34}+1)e^{(L_F+1)t}\|h\|_{H^{-1}(I)} , \; t\in (0,T]. \]
We end the proof with similar computations as in the proof of \cref{stability_of_test_function_space_by_the_semigroup}. 
$\square$
\begin{remark}
In fact, we have proved a slightly stronger result. Indeed, we have not used that $ \varphi\in \cE$ but only that $\varphi\in C^{1,0}_b(C(I)\times \R^k)$. Therefore
for $t>0$, $P_t$ maps $C^{1,0}_b(C(I)\times \R^k)$ into $\cE$.
\end{remark}

\begin{corollary}\label{characterization_of_the_full_generator_on_cE}\textbf{\texttt{(One more characterization of $\cA^\infty$)}}\vspace{0.1cm}\par 
\textbf{(i)} For all $\mu>0$, $\R(\mu-\cA^\infty)=\cE$, and \[ (\mu-\cA^\infty)^{-1}\varphi=\int_0^\infty e^{-\mu t}P_t\varphi dt \hspace{0.2cm}\forall \varphi\in\cE. \]
 
\textbf{(ii)} If $\mu>K=\max\{M,\tilde  M\}$, then for all $\varphi\in\cE$, $\psi:=(\mu-\cA^\infty)^{-1}\varphi\in\cE$.
\end{corollary}

\noindent \textit{Proof.} Since the domain of $\cA^\infty$ is $\cE$, under our assumption we already know that \[ \bar{\cA^\infty}\varphi(u)=\cA^\infty\varphi(u)=\hat{\cA^\infty}\varphi(u),\hspace{0.5cm}\forall\varphi\in\cE,\hspace{0.1cm}u\in E.\] Therefore, we conclude the part (i) by the Proposition 5.1, Section 5, Part 1 of \cite{Kurtz1986}, since the semigroup $(P_t)$ is measurable and of contraction on $\cE$ by $\cref{stability_of_test_function_space_by_the_semigroup}$, and observing that, thanks to Fubini theorem, the condition \[ P_s\int_0^\infty e^{-\mu t}P_t\varphi dt=\int_0^\infty e^{-\mu t}P_{s+t}\varphi dt \] holds for all $\varphi\in\cE$, $\mu>0$ and $s\geq0$.\par 

{Concerning part (ii), take $\varphi\in\cE$. From (i)}, it is clear that $\psi:=(\mu-\cA^\infty)^{-1}$ is bounded. Moreover, a derivation under the integral shows that $\psi$ is of class $C^1$ w.r.t. its first variable. At last, for all $\alpha\in B$ and $\nu\in K$, it follows from $\cref{stability_of_test_function_space_by_the_semigroup}$ that
\begin{align*}
|D_\alpha\psi(\alpha,\nu)\cdot h| & =\left|\int_0^\infty e^{-\mu t}D_\alpha P_t\psi(\alpha,\nu)\cdot hdt\right|\\
& \leq c\Vert \psi\Vert_\cE\left(\int_0^\infty e^{-(\mu-M)t}dt\right)\Vert h\Vert_\infty\\
& \leq c\Vert \psi\Vert_\cE\frac{1}{\mu-M}\Vert h\Vert_\infty,
\end{align*}
which yields boundedness for the differential of $\psi$. 
Similarly:
\begin{align*}
|D_\alpha\psi(\alpha,\nu)\cdot h| 
& \leq \tilde c\Vert \psi\Vert_\cE\left(\int_0^\infty (t^{-\frac34}+1)e^{-(\mu-\tilde M)t}dt\right)\Vert h\Vert_{H^{-1}(I)}\\
& \leq C\Vert \psi\Vert_\cE\frac{1}{(\mu-\tilde M)^{\frac14}}\Vert h\Vert_{H^{-1}(I)},
\end{align*}
for some constant $C$.
$\square$\\

We now use a classical argument to prove uniqueness. Let $\tilde{P}_{u}$ be another solution of the martingale problem for $\cA^\infty$. Let $\varphi\in\cE$, $\mu\in K$ and $\psi:=(\mu-\cA^\infty)^{-1}\varphi$. Then \[ \psi(v(t))-\psi(u)-\int_0^t\cA^\infty\psi(v(s))ds \] is a $\tilde{P}_u-$martingale. In particular, \[ \tilde{\E}_u\left(\psi(v(t)-\int_0^t\cA^\infty\psi(v(s))ds\right)=\psi(u). \] Multiply this identity by $\mu e^{-\mu t}$ and integrate on $[0,\infty)$ yield 
\[ \tilde{\E}_u\left(\int_0^\infty e^{-\mu t}\varphi(v(t))dt\right)=\psi(u)=\int_0^\infty e^{-\mu t}P_t\varphi(u)dt, \] for $\mu>K$, since $\cA^\infty$ is the full generator $\hat{\cA^\infty}$. By injectivity of the Laplace transform, this implies \[ \tilde{E}_u[\varphi(v(t))]=P_t\varphi(u)=\E_u[\varphi(v(t))], \] for almost all $t\geq0$. 

We have proved that the laws of the solutions to the martingale problem are the same at every fixed time $t$ in a dense set of $\R_+$. This implies uniqueness for the martingale problem (see \cite{Billingsley1999}, Section 14), and the proof of $\cref{well_posedness_of_the_martingale_problem_for_inf_dim_PDMP}$ ends. $\blacksquare$

\vspace{0.3cm}
\bibliographystyle{alpha}
\bibliography{MSSM_Hybrid_Limit_Supremum_Norm}

\begin{thebibliography}{CDMR12}

\bibitem[AD16]{DunErbZyg2016}
K.~Zygalakis A.~Duncan, R.~Erban.
\newblock Hybrid framework for the simulation of stochastic chemical kinetics.
\newblock {\em Journal of Computational Physics, Elsevier}, 326:398--419, 2016.

\bibitem[AT80]{Arnold1980}
L.~Arnold and M.~Theodosopulu.
\newblock Deterministic limit of the stochastic model of chemical reactions
  with diffusion.
\newblock {\em Adv. Appl. Prob.}, 12:367--379, 1980.

\bibitem[Bea09]{beale2009}
J.T. Beale.
\newblock Smoothing properties of implicit finite difference methods for a
  diffusion equation in maximum norm.
\newblock {\em SIAM Journal on Numerical Analysis}, 47(4), 2009.

\bibitem[Bil99]{Billingsley1999}
P. Billingsley.
\newblock {\em Convergence of Probability Measures, Second Edition}.
\newblock Wiley series in Probability and Statistics, 1999.

\bibitem[Blo87]{Blount1987}
D.~J. Blount.
\newblock {\em Comparison of a stochastic model of a chemical reaction with
  diffusion and the deterministic model}.
\newblock Ph.d. thesis, The University of Wisconsin-Madison, 1987.

\bibitem[Blo92]{blount1992}
D.~J. Blount.
\newblock Law of large numbers in the supremum norm for a chemical reaction
  with diffusion.
\newblock In {\em The Annals of Applied Probability}, volume~2, pages 131--141.
  1992.

\bibitem[BR11]{Riedler2011bis}
E.~Buckwar and M.~G. Riedler.
\newblock Exact modeling of neuronal membranes including spatio-temporal
  evolution.
\newblock {\em J. Math. Bio.}, 63((6)):1053--1091, 2011.

\bibitem[CDMR12]{Arnaud2012}
A.~Crudu, A.~Debussche, A.~Muller, and O.~Radulescu.
\newblock Convergence of stochastic gene networks to hybrid piecewise
  deterministic processes.
\newblock {\em The Anals of Applied Probability}, 22(5):1822--1859, 2012.

\bibitem[CDR09]{Arnaud2009}
A.~Crudu, A.~Debussche, and O.~Radulescu.
\newblock Hybrid stochastic simplifications for multiscale gene networks.
\newblock {\em BMC Systems Biology}, 3:89, Septembre 2009.

\bibitem[CH98]{Cazenave1998}
T. Cazenave and A. Haraux.
\newblock {\em An Introduction to Semilinear Evolution Equations}.
\newblock Clarendon Press - Oxford, 1998.

\bibitem[Dav84]{Davis1984}
M.~Davis.
\newblock Piecewise-deterministic markov processes: A general class of
  non-diffusion stochastic models.
\newblock {\em Journal of the Royal Statistical Society}, Series B
  (Methodological):353--388, 1984.

\bibitem[Dav93]{Davis1993}
M.~H.~A. Davis.
\newblock Markov models and optimization.
\newblock In Chapman and London Hall, editors, {\em Monographs on Statistics
  and Applied Probability}, volume~49. 1993.

\bibitem[DN17]{DebusscheNankep2017}
A.~{Debussche} and M.~J. {Nguepedja Nankep}.
\newblock {A Law of Large Numbers in the Supremum Norm for a Multiscale
  Stochastic Spatial Gene Network}.
\newblock {\em ArXiv e-prints}, November 2017.

\bibitem[EK86]{Kurtz1986}
S.~N. Ethier and T.~G. Kurtz.
\newblock {\em Markov Processes, Characterization and Convergence}.
\newblock John Wiley and Sons, Inc, 1986.

\bibitem[G{\'e}n13]{Genadot2013}
A. G{\'e}nadot.
\newblock {\em \'Etude multi-\'echelle de mod\`eles probabilistes pour les
  syst\`emes excitables avec composante spatiale}.
\newblock PhD thesis, Universit\'e Pierre et Marie Curie, 2013.

\bibitem[Hen81]{Henry1981}
D. Henry.
\newblock Geometric theory of semilinear parabolic equations.
\newblock In {\em Lecture Notes in Mathematics}. Springer, 1981.

\bibitem[JS87]{Jacod1987}
J. Jacod and A.~N. Shiryaev.
\newblock {\em Limit Theorems for Stochastic Processes}.
\newblock Springer-Verlag, Berlin Heidelberg GmbH, 1987.

\bibitem[Kat66]{Kato1966}
T.~Kato.
\newblock {\em Perturbation Theory for Linear Operators}.
\newblock Springer-Verlag, Berlin, 1966.

\bibitem[KK19]{Kouegou2019}
B. Kouegou~Kamen.
\newblock {\em Grandes d\'eviations dans des mod\`eles de biologie et des
  \'epid\'emies}.
\newblock Ph.d. thesis, Aix-Marseille Universit\'e, 2019.

\bibitem[Kot86]{Kotelenez1986}
P.~Kotelenez.
\newblock Law of large numbers and central limit theorem for linear chemical
  reactions with diffusion.
\newblock In {\em The Annals of Probability}, volume~14, pages 173--193.
  Universit\"at Bremen, 1986.

\bibitem[Kot88]{Kotelenez1988bis}
P.~Kotelenez.
\newblock A stochastic reaction-diffusion model.
\newblock In {\em University of Ultrecht}. 1988.

\bibitem[Kur70]{Kurtz1970}
T.~G. Kurtz.
\newblock Solutions of ordinary differential equations as limits of pure jump
  markov processes.
\newblock {\em J. Appl. Prob.}, 7:49--58, 1970.

\bibitem[Kur71]{Kurtz1971}
T.~G. Kurtz.
\newblock Limit theorems for sequences of jump markov processes approximating
  ordinary differential processes.
\newblock {\em J. Appl. Prob.}, 8:344--356, 1971.

\bibitem[NCS15]{NoelKarCheSch2015}
A. Noel, K.~C. Cheung, and R. Schober.
\newblock Multi-scale stochastic simulation for diffusive molecular
  communication.
\newblock In {\em Communication ICC}, London, UK, June 2015. IEEE.

\bibitem[NN18]{Nankep2018}
M.~J. Nguepedja~Nankep.
\newblock {\em Mod\'elisation stochastique de syst\`emes biologiques
  multi-\'echelles et inhomog\`enes en espace}.
\newblock Ph.d. thesis, Ecole Normale Sup\'erieure de Rennes, 2018.

\bibitem[NPY19]{NziPardouxYeo2019}
M.~{N'zi}, E.~{Pardoux}, and T.~{Yeo}.
\newblock {A SIR Model of a Refining Spatial Grid I: Law of Large Numbers}.
\newblock {\em Applied Mathematics and Optimization}, pages 1--37, 2019.

\bibitem[Rie11]{Riedler2011}
M.~G. Riedler.
\newblock {\em Spatio-temporal Stochastic Hybrid Models of Biological Excitable
  Membranes}.
\newblock PhD thesis, Heriot-Watt University, 2011.

\bibitem[RMC07]{RadMulCru2007}
O.~Radulescu, A.~Muller, and A.~Crudu.
\newblock Th{\'e}or{\`e}mes limites pour des processus de markov {\`a} sauts.
 Synth{\`e}se des r{\'e}sultats et applications en biologie mol{\'e}culaire.
\newblock {\em Tech. Sci. Inform.}, 26:443--469, 2007.

\bibitem[RTW12]{Riedler2012}
M.~G. Riedler, M.~Thieullen, and G.~Wainrib.
\newblock Limit theorems for infinite-dimensional piecewise deterministic
  markov processes. applications to stochastic excitable membrane models.
\newblock {\em Electron. J. Probab.}, 17(55):1--48, 2012.

\bibitem[Yeo19]{Yeo2019}
T. Yeo.
\newblock {\em Mod\`eles Stochastiques d'\'epid\'emies en Espace Discret et
  Continu : Loi des Grands Nombres et Fluctuations}.
\newblock Ph.d. thesis, Aix-Marseille Universit\'e, Universit\'e F\'elix
  Houphou\"et Boigny, 2019.

\end{thebibliography}

\end{document}